\documentclass[12pt]{article}
\usepackage{graphicx, amsmath}
\usepackage{amsfonts}
\usepackage{verbatim}
\usepackage{latexsym}
\makeatletter
\@addtoreset{equation}{section}

\begin{document}

\newcommand{\E}{\mathbb{E}}
\newcommand{\PP}{\mathbb{P}}
\newcommand{\RR}{\mathbb{R}}

\newtheorem{theorem}{Theorem}[section]
\newtheorem{lemma}{Lemma}[section]
\newtheorem{coro}{Corollary}[section]
\newtheorem{defn}{Definition}[section]
\newtheorem{assp}{Assumption}[section]
\newtheorem{expl}{Example}[section]
\newtheorem{prop}{Proposition}[section]
\newtheorem{remark}{Remark}[section]

\newcommand\tq{{\scriptstyle{3\over 4 }\scriptstyle}}
\newcommand\qua{{\scriptstyle{1\over 4 }\scriptstyle}}
\newcommand\hf{{\textstyle{1\over 2 }\displaystyle}}
\newcommand\hhf{{\scriptstyle{1\over 2 }\scriptstyle}}

\newcommand{\eproof}{\indent\vrule height6pt width4pt depth1pt\hfil\par\medbreak}

\def\a{\alpha} \def\g{\gamma}
\def\e{\varepsilon} \def\z{\zeta} \def\y{\eta} \def\o{\theta}
\def\vo{\vartheta} \def\k{\kappa} \def\l{\lambda} \def\m{\mu} \def\n{\nu}
\def\x{\xi}  \def\r{\rho} \def\s{\sigma}
\def\p{\phi} \def\f{\varphi}   \def\w{\omega}
\def\q{\surd} \def\i{\bot} \def\h{\forall} \def\j{\emptyset}
\def\R{\mathbb R} 
\def\be{\beta} \def\de{\delta} \def\up{\upsilon} \def\eq{\equiv}
\def\ve{\vee} \def\we{\wedge}
\def\dd{\delta}
\def\F{{\cal F}}
\def\T{\tau} \def\G{\Gamma}  \def\D{\Delta} \def\O{\Theta} \def\L{\Lambda}
\def\X{\Xi} \def\Si{\Sigma} \def\W{\Omega}
\def\M{\partial} \def\N{\nabla} \def\Ex{\exists} \def\K{\times}
\def\V{\bigvee} \def\U{\bigwedge}

\def\1{\oslash} \def\2{\oplus} \def\3{\otimes} \def\4{\ominus}
\def\5{\circ} \def\6{\odot} \def\7{\backslash} \def\8{\infty}
\def\9{\bigcap} \def\0{\bigcup} \def\+{\pm} \def\-{\mp}
\def\la{\langle} \def\ra{\rangle}

\def\tl{\tilde}
\def\trace{\hbox{\rm trace}}
\def\diag{\hbox{\rm diag}}
\def\for{\quad\hbox{for }}
\def\refer{\hangindent=0.3in\hangafter=1}

\newcommand\wD{\widehat{\D}}
\newcommand{\ka}{\kappa_{10}}

\title{
\bf {Tamed EM scheme of Neutral Stochastic Differential Delay Equations}
}
\author
{
{\bf Yanting Ji} \
and \ {\bf Chenggui Yuan}\thanks{Contact e-mail address: C.Yuan@swansea.ac.uk}\\[1ex]
 Department of Mathematics\\
Swansea University, Swansea,  SA2 8PP, U.K.\\
}

\date{}

\maketitle
\begin{abstract}
In this paper, we investigate the convergence of the tamed Euler-Maruyama (EM) scheme for  a class of neutral stochastic differential delay equations. The strong convergence results of  the tamed EM scheme are presented under global and local non-Lipschitz conditions, respectively. Moreover, under a global Lipschitz condition, we provide the rate of the  convergence of tamed EM, which is smaller than the rate convergence of classical EM scheme  one half. 
\\

{\it MSC 2010\/}: 65C30 (65L20, 60H10) 
	\smallskip
	
{\it Key Words and Phrases}: neutral stochastic differential delay equations, non-Lipschitz, monotonicity, tamed EM scheme, rate of convergence, pure jumps, Poisson processes.
	
\end{abstract}
\section{Introduction}

The Euler-Maruyama scheme is of vital importance in numerical approximation for stochastic differential equations (SDEs). In \cite{Klo}, Kloeden and Platen illustrated that, if the coefficients of an SDE are globally Lipschitz continuous, then the EM approximation converges to the exact solution of the SDE in both strong and weak sense, the convergence rates for both cases are provided as well. In the same book, they also mentioned that the Milstein scheme converges to the exact solution of SDE in both strong and weak sense with different orders under certain conditions including the global Lipschitz condition. It is the first time that  Higham, Mao and Stuart  \cite{Hig} established
strong convergence results under the super-linear condition and the moment boundedness condition, however, it remained an open question whether the moment of the EM approximation is bounded within finite time if the coefficients of an SDE are not globally Lipschitz
continuous.
Recently,  Hutzenthaler, Jentzen and Kloeden  \cite{Hut} have found that once the global Lipschitz condition was replaced by the super-linear condition, the moment of the EM scheme could be infinity within the finite time. To tackle this problem, in the paper \cite{Jen}, Hutzenthaler, Jentzen and Kloeden introduced a new approximation scheme, which is the so-call tamed EM scheme. By employing the tamed scheme, the drift coefficient is tamed so that it is uniformly bounded. With such an approach, it is proved that the tamed EM scheme converges to the exact solution of the SDE under the super-linear condition of the drift coefficients. The tamed EM scheme is later extended to SDEs with locally Lipschitz condition of the diffusion by Sabanis \cite{Sab}.

On the other hand, stochastic differential delay equations (SDDEs) and neutral stochastic differential delay equations (NSDDEs) describe a wide variety of natural and man-made
systems. For the theories and applications of SDDEs and NSDDEs, we here only mention \cite{Arr,BY,Bell,Bou,Els,Es,HMS1,jwy09,Liao,MYZ,McW,Moh}, to name a few.
Since most SDDEs and NSDDEs can not be
solved explicitly, numerical methods have become essential.
Recently, an extensive literature has emerged in investigating the strong
convergence, weak convergence or sample path convergence  of
numerical schemes for SDDEs and NSDDEs,
for example, \cite{bbmy11,hu96,Hu, jwy09,KP00,ms03}.
 We should  point out that the strong
convergence of EM schemes for SDDEs is, in general, discussed under
a linear growth condition or bounded moments of analytic and
numerical solutions, e.g., \cite{jwy09,KP00,ms03}. However, similar to the SDEs case, it remained an open question whether the EM scheme converges to the exact solution if the coefficients of the SDDEs and NSDDEs are under the super-linear condition. The main aim of this paper is to answer this question by extending the tamed EM method to NSDDEs. Because of the neutral term, the technical details increased significantly.

The rest of this paper will be organised as follows. In section 2, some notation and preliminaries are introduced. In section 3, $p$-th moment boundedness, convergence of EM scheme under global and local monotonicity conditions are proven, moreover, the rate of convergence is provided under a global conditions.  In section 4, we present similar results  as in section 3 while the Brownian motion is replaced by the pure jump processes.

\section{Preliminaries}
Throughout this paper, let $(\Omega,\mathcal{F},\PP)$ be a complete probability space with a filtration $\{\mathcal{F}_t\}_{t\geq 0}$ satisfying the usual condition (i.e. it is right continuous and $\mathcal{F}_{0}$ contains all $\PP$-null sets).
Let $\tau>0$ be a constant and denote $C([-\tau,0];\mathbb{R}^{n})$ the space of all continuous functions from $[-\tau,0]$ to $\mathbb{R}^{n}$ with the norm $\|\phi\|=\sup_{-\tau\leq\theta\leq 0}|\phi(\theta)|$. 
 Let $B(t)$ be a standard $m$-dimensional Brownian motion.

Consider an $n$-dimensional neutral stochastic differential delay equation
\begin{equation}\label{NSDDE}
d[X(t)-D(X(t-\tau))]  = b(X(t),X(t-\tau))dt+\sigma(X(t),X(t-\tau))dB(t),
\end{equation}
on $t \geq 0$, where
\begin{equation*}
 D:\mathbb{R}^{n} \to \mathbb{R}^{n},\,
 b:\mathbb{R}^{n}\times \mathbb{R}^{n}\to \mathbb{R}^{n},\,
\sigma:\mathbb{R}^{n}\times \mathbb{R}^{n} \to \mathbb{R}^{n\times m}, 
\end{equation*}
we assume that $D$, $b$ and $\sigma$ are Borel-measurable  and the initial data satisfies the following condition: for any $p\ge 2$
\begin{equation}\label{initial}
\{X(\theta):-\tau\leq \theta\leq 0\} = \xi \in L^p_{\mathcal{F}_{0}}([-\tau,0];\mathbb{R}^{n}),
\end{equation}
that is  $\xi $ is an $\mathcal{F}_{0}$-measurable $C([-\tau,0];\mathbb{R}^{n})$-valued random variable and  $E\|\xi\|^p <\infty.$
Now, fix $T >\tau > 0$, without loss of generality, we assume that $T$ and $ \tau$ are rational numbers, and the step size $h\in (0, 1)$ be fraction of $\tau$ and $T,$  so that there exist two positive integers $M, \bar{M}$ such that $h=T/M=\tau/\bar{M}.$ 
Throughout the paper, unless otherwise stated, $C >0$ denotes a generic constant, independent of $h$, whose value may change from line to line.

For the future use,  we assume that:

\begin{description}
\item[(A1)]There exists a positive constant $\tilde{K}$ such that
\begin{equation}\label{WeakCorecivity}
\langle x-D(y),b(x,y)\rangle \vee ||\sigma(x,y)||^{2}\leq \tilde{K}(1+|x|^{2}+|y|^2),
\end{equation}
for $\forall$ $x,y \in \mathbb{R}.$\\
\item[(A2)]  $D(0)=0$ and there exists a constant $\kappa \in (0,1)$ such that
\begin{eqnarray}\label{NeutralBound}
	|D(x)-D(\bar{x})| \leq \kappa|x-\bar{x}| \mbox{ for all }x,y \in \mathbb{R}^n.
\end{eqnarray}
\item[(A3)] For any $R>0,$ there exist two  positive constants $\tilde{K}_{R}$ and $K_{R}$ such that
\begin{equation}\label{LocalLipschitzII}
\begin{split}
& \langle x-D(y)-\bar{x}+D(\bar{y}), b(x,y)-b(\bar{x},\bar{y})\rangle\vee||\sigma(x,y)-\sigma(\bar{x},\bar{y})||^2 \\
&\leq \tilde{K}_{R} (|x-\bar{x}|^{2}+|y-\bar{y}|^2) ,
\end{split}
\end{equation}
for all $|x|\vee|y|\vee|\bar{x}|\vee|\bar{y}| \leq R,$\\
and 
\begin{equation}
\sup_{|x|\vee |y|\leq R}|b(x,y)| \leq K_{R}.
\end{equation}

\item[(A4)]There exist two positive constants $l$ and $L$ such that,
\begin{equation*}
\langle x-D(y)-\bar{x}+D(\bar{y}), b(x,y)-b(\bar{x},\bar{y})\rangle+ ||\sigma(x,y)-\sigma(\bar{x},\bar{y})||^2 \leq L(|x-\bar{x}|^2 + |y-\bar{y}|^2),
\end{equation*}
and
\begin{equation*}
|b(x,y)-b(\bar{x},\bar{y})|\leq L(1+|x|^{l}+|y|^{l}+|\bar{x}|^{l}+|\bar{y}|^{l})(|x-\bar{x}|+|y-\bar{y}|),
\end{equation*}
for all $x,$ $y,$ $\bar{x}$ and $\bar{y} \in \mathbb{R}^{n}.$
\item[(A5)] For any $s,t\in[-\tau,0]$ and $q>0,$ let
$\E\|\xi(t)-\xi(s)\|^{q}\leq L|t-s|^q.$ 
\end{description}
\smallskip

\begin{remark}\label{EUTheorem}
{\rm Assume that {\bf (A1)}-{\bf (A3)} hold, then NSDDE \eqref{NSDDE} with initial \eqref{initial} admits a unique strong global solution $X(t),$ $t\in[0,T].$ The proof details of such existence and uniqueness result can be found in \cite{Ji}. If assumption {\bf (A3)} is replaced by {\bf (A4)}, the theorem of  existence and uniqueness  still holds. }
\end{remark}

We now define
\begin{eqnarray}\label{SpecialTamedcoefficients}
b_{h}(x,y):= \frac{b(x,y)}{1+h^{\alpha}|b(x,y)|},
\end{eqnarray}
for all  $x,y\in \mathbb{R}^{n}$ and $\alpha \in (0,\frac{1}{2}].$  The reason why such value of $\alpha$ is chosen  will be revealed in Section 3.
By observation, one has
\begin{equation*}
|b_{h}(x,y)| \leq \min(h^{-\alpha},|b(x,y)|).
\end{equation*}
\begin{remark}\label{TamedSol}
{\rm It is easy to verify that $b_{h}(x,y) = \frac{b(x,y)}{1+h^{\alpha}|b(x,y)|}$ satisfies all assumptions {\bf (A1)}, {\bf (A3)} and {\bf (A4)}.
For the assumption {\bf (A1)}, we have
\begin{equation*}
\begin{split}
2\langle x-D(y),b_{h}(x,y)\rangle & = \frac{2}{1+h^{\alpha}|b(x,y)|}\langle x-D(y),b(x,y)\rangle\\
&\leq \frac{\tilde{K}}{1+h^{\alpha}|b(x,y)|}(1+|x|^{2}+|y|^2)\\
&\le \tilde{K}(1+|x|^{2}+|y|^2).
\end{split}
\end{equation*}
For the assumption {\bf (A3)}, which is a local assumption, we can derive that
\begin{equation*}
\begin{split}
& 2\langle x-D(y)-\bar{x}+D(\bar{y}), b_{h}(x,y)-b_{h}(\bar{x},\bar{y})\rangle \\
&\leq \bar{K}_{R} (|x-\bar{x}|^{2}+|y-\bar{y}|^2) ,
\end{split}
\end{equation*}
where $\bar{K}_{R} = \frac{\tilde{K}_{R}}{1+h^{\alpha}(|b(\bar{x},\bar{y})|\wedge|b(x,y)|)}.$ The verification detail for assumption {\bf (A4)} is omitted, since it is similar to {\bf (A3)}.}
\end{remark}

Now, we can define the discrete-time tamed EM scheme:\\
For every integer $n = -\bar{M},\cdots,0,$ $Y_{h}^{(n)} = \xi(nh).$ For every integer $n = 1, \cdots, M-1,$
\begin{equation}\label{DEulerschemeStep}
Y_{h}^{(n+1)}-D(Y_{h}^{(n+1-\bar{M})})= Y_{h}^{(n)}-D(Y_{h}^{(n-\bar{M})})+b_{h}(Y_{h}^{(n)},Y_{h}^{(n-\bar{M})})h+\sigma(Y_{h}^{(n)},Y_{h}^{(n-\bar{M})})\Delta B^{(n)}_{h},
\end{equation}
where $\Delta B^{(n)}_{h} = B((n+1)h)-B(nh).$ 
The discrete-time tamed EM scheme \eqref{DEulerschemeStep} can be rewritten as
\begin{equation}\label{DEulerschemeSum}
\begin{split}
Y_{h}^{(n+1)}=&D(Y_{h}^{(n+1-\bar{M})})+ \xi(0)-D(\xi(-\tau))+\sum_{i=1}^{n}b_{h}(Y_{h}^{(i)},Y_{h}^{(i-\bar{M})})h\\
&+\sum_{i=1}^{n}\sigma(Y_{h}^{(i)},Y_{h}^{(i-\bar{M})})\Delta B^{(i)}_{h}.
\end{split}
\end{equation}
For $t\in[nh,(n+1)h),$  we set  $\bar{Y}(t): = Y_{h}^{(n)}.$  Since $\tau=\bar M h ,$ $\bar{Y}(t-\tau) = Y_{h}^{(n-\bar{M})}.$  
For the sake of simplicity, we define the corresponding continuous-time tamed EM approximate solution $Y(t)$ as follows: 
For any $\theta\in[-\tau,0],$ $Y(\theta) = \xi(\theta).$
For any $t \in [0,T],$
\begin{equation}\label{CEulerscheme}
\begin{split}
Y(t)=& D(\bar{Y}(t-\tau))+\xi(0)-D(\xi(-\tau))+\int_{0}^{t}b_{h}(\bar{Y}(s),\bar{Y}(s-\tau))ds\\
&+\int_{0}^{t}\sigma(\bar{Y}(s),\bar{Y}(s-\tau))dB(s) .
\end{split}
\end{equation}
Noting that for any $t\in[0,T],$ there exists a positive integer $n,$ $0\leq n\leq M-1,$ such that $nh\leq t< (n+1)h,$ we have
\begin{equation}\label{CEulerschemeModification}
\begin{split}
Y(t) &=D(\bar{Y}(t-\tau))+\xi(0)-D(\xi(-\tau))+\int_{0}^{nh}b_{h}(\bar{Y}(s),\bar{Y}(s-\tau))ds\\
& +\int_{0}^{nh}\sigma(\bar{Y}(s),\bar{Y}(s-\tau))dB(s)+\int_{nh}^{t}b_{h}(\bar{Y}(s),\bar{Y}(s-\tau))ds\\
& +\int_{nh}^{t}\sigma(\bar{Y}(s),\bar{Y}(s-\tau))dB(s)\\
&= Y(nh)+\int_{nh}^{t}b_{h}(\bar{Y}(s),\bar{Y}(s-\tau))ds+\int_{nh}^{t}\sigma(\bar{Y}(s),\bar{Y}(s-\tau))dB(s).\\
\end{split}
\end{equation}
This means  the continuous-time tamed EM approximate solution $Y(t)$ coincides with the discrete-time tamed approximation solution $\bar{Y}(t)$ at grid points $t = nh,$ $n = 0,1,\cdots,M-1.$ 

\section{Tamed EM Method of NSDDEs driven by  Brownian Motion}

In this section, we show that the tamed EM scheme converges to the exact solution under certain conditions,  i.e. we have the following main results: 
\begin{theorem}\label{ConvergenceRate}
Suppose that {\bf(A1)}, {\bf(A2)}, {\bf(A4)} and {\bf(A5)} hold, then the tamed EM scheme $(\ref{CEulerscheme})$ converges to the exact solution of $(\ref{NSDDE})$ such that for any $p\geq 2,$
\begin{equation}\label{MainResultConvergenceRate}
E\Big[\sup_{0 \leq t \leq T}|X(t)-Y(t)|^{p}\Big] \leq Ch^{\alpha p}.
\end{equation}
\end{theorem}
The next theorem states that the convergence result still holds if the global monotonic condition  {\bf(A4)} is replaced by its local counterpart {\bf(A3)}. However, we are unable to provide the convergence rate under this weaker condition:
\begin{theorem}\label{Convergence}
Suppose {\bf(A1)}-{\bf(A3)} and {\bf(A5)} hold, then the tamed EM scheme $(\ref{CEulerscheme})$ converges to the exact solution of $(\ref{NSDDE})$ such that for any $p\geq2,$
\begin{equation}\label{MainResultConvergence}
\lim_{h \to 0}\E\Big[\sup_{0 \leq t \leq T}|X(t)-Y(t)|^{p}\Big] = 0.
\end{equation}
\end{theorem}

\subsection{Moment Bounds}
In order to prove our main results, in this subsection we investigate the boundedness of moments of  both exact solution and EM approximation. 
\begin{lemma}\label{Lemma1}
Consider the continuous-time tamed EM scheme given by equation $(\ref{CEulerschemeModification}).$ If for some $p \geq 2,$
\begin{equation}\label{AssumptionLemma1}
\sup_{0\leq t\leq T}\E(|Y(t)|^{p}) \le C,
\end{equation}
and {\bf (A1)} hold, then
 it holds that 
\begin{equation}\label{EstimateIntervalError}
\E\bigg[\sup_{0\le n \le M-1}\sup_{nh\leq t \leq (n+1)h}|Y(t)-Y(nh)|^p\bigg]\leq Ch^{p/2},
\end{equation}
and
\begin{equation}\label{ErrorWithDrift}
\E\bigg[\sup_{0\le n \le M-1}\sup_{nh\leq t \leq (n+1)h}|Y(t)-Y(nh)|^p|b_{h}(\bar{Y}(t),\bar{Y}(t-\tau))|^{p}\bigg] \leq C.
\end{equation}
\end{lemma}
\textbf{Proof}: 
By the definition of the tamed EM scheme, we have for $nh\le t <(n+1)h,$
\begin{equation*}
\begin{split}
\E\bigg[\sup_{nh\leq t \leq (n+1)h}|Y(t)-Y(nh)|^p\bigg] &= \E\bigg[\sup_{nh\leq t \leq (n+1)h}\bigg|\int_{nh}^{t}b_{h}(\bar{Y}(s),\bar{Y}(s-\tau))ds\\
&~~~~~~+\int_{nh}^{t}\sigma(\bar{Y}(s),\bar{Y}(s-\tau))dB(s)\bigg|^{p}\bigg].
\end{split}
\end{equation*}
 Therefore, due to H\"older's inequality,
\begin{equation}\label{Lemma1Observation}
\begin{split}
\E\bigg[\sup_{nh\leq t \leq (n+1)h}|Y(t)-Y(nh)|^p\bigg]&\leq 2^{p-1} h^{p-1}\E\bigg[\int_{nh}^{(n+1)h}\bigg|b_{h}(\bar{Y}(s),\bar{Y}(s-\tau))\bigg|^{p}ds\bigg]\\
&+2^{p-1}\E\bigg[\sup_{nh\leq t \leq (n+1)h}\bigg|\int_{nh}^{t}\sigma(\bar{Y}(s),\bar{Y}(s-\tau))dB(s)\bigg|^{p}\bigg].
\end{split}
\end{equation}
Using the Burkholder-Davis-Gundy(BDG) inequality \cite[Theorem 1.7.3 page 40]{Mao} and  $(\ref{AssumptionLemma1}),$  for some $p \geq 2,$ we can derive that for $nh\le t <(n+1)h,$
\begin{equation*}
\begin{split}
&\E\bigg[\sup_{nh\leq t \leq (n+1)h}\big|\int_{nh}^{t}\sigma(\bar{Y}(s),\bar{Y}(s-\tau))dB(s)\big|^{p}\bigg]\\
&\leq C\E\left[\int_{nh}^{(n+1)h}\|\sigma(\bar{Y}(s),\bar{Y}(s-\tau))\|^{2}ds\right]^{p/2}\\
&\leq C\E\left[\int_{nh}^{(n+1)h}(1+|\bar{Y}(s)|^{2}+|\bar{Y}(s-\tau)|^{2})ds\right]^{p/2}\\
&\leq Ch^{p/2}.
\end{split}
\end{equation*}
This, together with  $|b_{h}(\bar{Y}(s),\bar{Y}(s-\tau))|\leq h^{-\alpha},$ yields
\begin{equation}\label{ay}
\E\bigg[\sup_{nh\leq t \leq (n+1)h}|Y(t)-Y(nh)|^p\bigg] \leq  2^{p-1}h^{(1-\alpha)p}+Ch^{p/2} \leq Ch^{p/2}.
\end{equation}
Therefore $(\ref{EstimateIntervalError})$ holds.
Moreover,
\begin{equation}\label{ay1}
\begin{split}
&\E\bigg[\sup_{nh\leq t \leq (n+1)h}|Y(t)-Y(nh)|^p|b_{h}(\bar{Y}(s),\bar{Y}(s-\tau))|^{p}\bigg] \\
&\leq \E\bigg[\sup_{nh\leq t \leq (n+1)h}|Y(t)-Y(nh)|^p\bigg]h^{-\alpha p} \\
&\leq Ch^{(1/2-\alpha)p}\leq C.
\end{split}
\end{equation}
In \eqref{ay} and \eqref{ay1}, we have used the fact that $\alpha\in (0,1/2].$ The proof is therefore complete.
 \hfill $\Box$

\begin{lemma}\label{yl3}
Assume that {\bf (A1)}, {\bf (A2)} and {\bf (A5)} hold. Then there exists a positive constant $C$ independent of $h$ such that for any $p\geq 2,$	
\begin{equation}\label{UpperBuondLemma3}
\E[\sup_{0\leq t \leq T}|X(t)|^{p}]\vee \E[\sup_{0\leq t \leq T}|Y(t)|^{p}]\leq C.
\end{equation}	
\end{lemma}
{\bf Proof.}\quad
By \cite[Lemma 4.4, p212]{Mao} and assumption {\bf(A2)}, we know that for any $p\geq 2,$
\begin{equation}\label{Lemma301}
\begin{split}
\sup_{0\leq t \leq T}|X(t)|^{p} \leq \frac{\kappa}{1-\kappa}||\xi||^{p}+\frac{1}{(1-\kappa)^p}\sup_{0\leq t \leq T}|X(t)-D(X(t-\tau))|^{p},
\end{split}
\end{equation}
An application of It\^o's formula yields
\begin{equation}\label{eq50}
\begin{split}
&|X(t)-D(X(t-\tau))|^{p}\leq |\xi(0)-D(\xi(-\tau))|^{p}\\
&+p\int_{0}^{t}|X(t)-D(X(t-\tau))|^{p-2}\langle X(s)-D(X(s-\tau)),b(X(s),X(s-\tau))\rangle ds\\
&+\frac{p(p-1)}{2}\int_{0}^{t}|X(t)-D(X(t-\tau))|^{p-2}||\sigma(X(s),X(s-\tau))||^{2}ds\\
&+p\int_{0}^{t}|X(t)-D(X(t-\tau))|^{p-2}(X(s)-D(X(s-\tau)))^{T}(\sigma(X(s),X(s-\tau))dB(s))\\
&=: |\xi(0)-D(\xi(-\tau))|^{p}+ H_{1}(t)+H_{2}(t)+H_{3}(t).
\end{split}
\end{equation}
Under {\bf(A1)} and {\bf(A2)}, one has
\begin{equation}\label{H1}
\begin{split}
&\E(\sup_{0\leq t\leq T}H_{1}(t))+\E(\sup_{0\leq t\leq T}H_{2}(t))\\
&\leq C\E\int_{0}^{T}|X(s)-D(X(s-\tau))|^{p-2}(1+|X(s)|^{2}+|X(s-\tau)|^{2})ds\\
&\leq C\E\int_{0}^{T}(|X(s)|^{p-2}+|D(X(s-\tau))|^{p-2})(1+|X(s)|^{2}+|X(s-\tau)|^{2})ds\\
&\leq C\E\int_{0}^{T}(1+|X(s)|^{p}+|X(s-\tau)|^{p})ds\leq C+C\int_{0}^{T}\E(\sup_{0\le u\le s}|X(u)|^p)ds.
\end{split}
\end{equation}
By the Burkholder-Davis-Gundy(BDG) inequality and the Young inequality, we derive that
\begin{equation}\label{H2}
\begin{split}
&\E(\sup_{0\leq t\leq T}H_{3}(t))\leq C\E\bigg(\int_{0}^{T}\big(|X(s)-D(X(s-\tau))|^{2p-2}\|\sigma(X(s),X(s-\tau))\|^{2}\big)ds\bigg)^{1/2}\\
&\leq C\E\bigg(\sup_{0\leq t\leq T}|X(t)-D(X(t-\tau))|^{p-1}\big(\int_{0}^{T}\|\sigma(X(s),X(s-\tau))\|^{2}ds\big)^{1/2}\bigg)\\
&\leq\frac{1}{4}\E(\sup_{0\leq t \leq T}|X(t)-D(X(t-\tau))|^{p})+C(\int_{0}^{T}\|\sigma(X(s),X(s-\tau))\|^{2}ds)^{p/2}\\
&\leq\frac{1}{4}\E(\sup_{0\leq t \leq T}|X(t)-D(X(t-\tau))|^{p})+C(\int_{0}^{T}(1+|X(s)|^{2}+|X(s-\tau)|^{2}ds)^{p/2}\\
&\leq \frac{1}{4}\E(\sup_{0\leq t \leq T}|X(t)-D(X(t-\tau))|^{p})+C\int_{0}^{T}\E(1+|X(s)|^{p}+|X(s-\tau)|^{p})ds\\
&\leq\frac{1}{4}\E(\sup_{0\leq t \leq T}|X(t)-D(X(t-\tau))|^{p})+C+C\int_{0}^{T}\E(\sup_{0\le u\le s}|X(u)|^p)ds.\\
\end{split}
\end{equation}
Now, substituting \eqref{H1} and \eqref{H2} into \eqref{eq50},  we then  have
\begin{equation}\label{eq51}
\E(\sup_{0\leq t\leq T}|X(t)-D(X(t-\tau))|^{p})\leq C+C\int_{0}^{T}\E(\sup_{0\le u\le s}|X(u)|^p)ds.
\end{equation}
This, together with \eqref{Lemma301}, implies that
\begin{equation*}
\E(\sup_{0\leq t\leq T}|X(t)|^p)\leq C+C\int_{0}^{T}\E(\sup_{0\le u\le s}|X(u)|^p)ds.
\end{equation*}
The desired assertion for the exact solution follows an application of Gronwall's inequality.

In order to estimate the $p$-th moment of the tamed EM scheme \eqref{CEulerscheme}, an inductive argument is used below. Firstly, we claim that there exists a constant $C$ such that:
\begin{equation}\label{y319}
\sup_{0\le t\le T}\E |Y(t)|^2 \le C.
\end{equation}
An application of the It\^o formula yields, 
\begin{equation}\label{H3001}
\begin{split}
&|Y(t)-D(\bar{Y}(t-\tau))|^{2} \\
&= |\xi(0)-D(\xi(-\tau))|^{2}+2\int_{0}^{t}\langle Y(s)-D(\bar{Y}(s-\tau)),b_{h}(\bar{Y}(s),\bar{Y}(s-\tau))\rangle ds\\
&+\int_{0}^{t}||\sigma(\bar{Y}(s),\bar{Y}(s-\tau))||^{2}ds+2\int_{0}^{t}\langle Y(s)-D(\bar{Y}(s-\tau)),\sigma(\bar{Y}(s),\bar{Y}(s-\tau)) dB(s)\rangle\\
&= |\xi(0)-D(\xi(-\tau))|^{2}+2\int_{0}^{t}\langle\bar{Y}(s)-D(\bar{Y}(s-\tau)),b_{h}(\bar{Y}(s),\bar{Y}(s-\tau))\rangle ds\\
&+2\int_{t_0}^{t}\langle Y(s)-\bar{Y}(s),b_{h}(\bar{Y}(s),\bar{Y}(s-\tau))\rangle ds+\int_{0}^{t}||\sigma(\bar{Y}(s),\bar{Y}(s-\tau))||^{2}ds\\
&+2\int_{0}^{t}\langle Y(s)-D(\bar{Y}(s-\tau)),\sigma(\bar{Y}(s),\bar{Y}(s-\tau))dB(s)\rangle \\
&=: |\xi(0)-D(\xi(-\tau))|^{2}+\hat{H}_{1}(t)+\hat{H}_{2}(t)+\hat{H}_{3}(t)+\hat{H}_{4}(t).\\
\end{split}
\end{equation}
By  {\bf(A1)} and {\bf(A2)},  we compute
\begin{equation}\label{H100}
\begin{split}
&\E(\hat H_{1}(t))+\E(\hat H_{3}(t))
\leq C\E\int_{0}^{t}1+|\bar Y(s)|^{2}+|\bar{Y}(s-\tau)|^{2}ds\\
&\leq C+C\int_{0}^{t}\sup_{0\le u\le s}\E(|Y(u)|^2)ds.
\end{split}
\end{equation}
Using the definition of $Y(s)$ and $\bar Y(s)$, we have
\begin{equation}\label{H300}
\begin{split}
&\E(\hat H_{2}(t))= 2\E\int_{0}^{t}\langle \int_{[\frac{s}{h}]h}^{s}b_{h}(\bar{Y}(r),\bar{Y}(r-\tau))dr, b_{h}(\bar{Y}(s),\bar{Y}(s-\tau))\rangle ds\\
&+2\E\int_{0}^{t}\langle\int_{[\frac{s}{h}]h}^{s}\sigma(\bar{Y}(r),\bar{Y}(r-\tau))dB(r),b_{h}(\bar{Y}(s),\bar{Y}(s-\tau))\rangle ds\\
& = 2\E\int_{0}^{t}\langle \int_{[\frac{s}{h}]h}^{s}b_{h}(\bar{Y}(r),\bar{Y}(r-\tau))dr, b_{h}(\bar{Y}(s),\bar{Y}(s-\tau))\rangle ds\\
&+2\E\int_{0}^{t}\langle\int_{[\frac{s}{h}]h}^{s}\sigma(\bar{Y}(r),\bar{Y}(r-\tau))dB(r)|{{\cal{F}}_{[\frac{s}{h}]h}},b_{h}(\bar{Y}(s),\bar{Y}(s-\tau))\rangle ds\\&= 2\E\int_{0}^{t}\langle \int_{[\frac{s}{h}]h}^{s}b_{h}(\bar{Y}(r),\bar{Y}(r-\tau))dr, b_{h}(\bar{Y}(s),\bar{Y}(s-\tau))\rangle ds\\&\leq Cth^{1-2\alpha}\leq C,
\end{split}
\end{equation} 
here again we have used $\alpha \in (0, 1/2].$ 
Using theYoung inequality, we have 
\begin{equation}\label{H304}
\begin{split}
&\sup_{0\leq t\leq T}\E(|Y(t)|^2)= \sup_{0\leq t\leq T}\E(|Y(t)-D(\bar{Y}(t-\tau))+D(\bar{Y}(t-\tau))|^2)\\
&\leq C+\frac{1}{2}\sup_{0\leq t\leq T}\E(|Y(t)|^2)+C\sup_{0\leq t\leq T}\E(|Y(t)-D(\bar{Y}(t-\tau))|^2),
\end{split}
\end{equation}
 Now, substituting \eqref{H100} and \eqref{H300}   into \eqref{H304}, we can derive that
\begin{equation*}
\begin{split}
&\sup_{0\leq t \leq T}\E(|Y(t)|^{2}) \le C+ C\big(\sup_{0\leq t \leq T}\E(|Y(t)-D(\bar{Y}(t-\tau))|^{2})\big)\\
&\leq C+C\int_0^T\sup_{0\le u\le s}\E(|Y(u)|^2)ds.
\end{split}
\end{equation*}
The assertion \eqref{y319} follows  an application of the Gronwall inequality.

In the sequel, we are now going to show that there exists a constant $C>0$ such that for any $p\geq 2,$
\begin{equation}\label{y320}
\E[\sup_{0 \le t \le T}|Y(t)|^p]\le C.
\end{equation}
Letting  $p=4 $ and  using the  It\^o formula, we have
\begin{equation}\label{EstimatedIntervalIto}
\begin{split}
&|Y(t)-D(\bar{Y}(t-\tau))|^{p}\leq |\xi(0)-D(\xi(-\tau))|^{p}+p\int_{0}^{t}|Y(s)-D(\bar{Y}(s-\tau))|^{p-2}\\
&~~~~~~~~~~~~\K\langle\bar{Y}(s)-D(\bar{Y}(s-\tau)),b_{h}(\bar{Y}(s),\bar{Y}(s-\tau))\rangle ds\\
&+p\int_{0}^{t}|Y(s)-D(\bar{Y}(s-\tau))|^{p-2}\langle Y(s)-\bar{Y}(s),b_{h}(\bar{Y}(s),\bar{Y}(s-\tau))\rangle ds\\
&+\frac{p(p-1)}{2}\int_{0}^{t}|Y(s)-D(\bar{Y}(s-\tau))|^{p-2}\|\sigma(\bar{Y}(s),\bar{Y}(s-\tau))\|^{2}ds\\
&+p\int_{0}^{t}|Y(s)-D(\bar{Y}(s-\tau))|^{p-2}\langle Y(s)-D(\bar{Y}(s-\tau)),\sigma(\bar{Y}(s),\bar{Y}(s-\tau))dB(s)\rangle \\
&=: |\xi(0)-D(\xi(-\tau))|^{p}+ \bar H_{1}(t)+\bar H_{2}(t)+\bar H_{3}(t)+\bar H_{4}(t).\\
\end{split}
\end{equation}
Under assumptions {\bf(A1)} and {\bf(A2)}, there exists a constant $C>0$ such that 
\begin{equation}\label{H101}
\begin{split}
&\E(\sup_{0\leq t\leq T}\bar H_{1}(t))+\E(\sup_{0\leq t\leq T}\bar H_{3}(t))\\
&\leq C\E\int_{0}^{T}|Y(s)-D(\bar{Y}(s-\tau))|^{p-2}(1+\bar |Y(s)|^{2}+|\bar{Y}(s-\tau)|^{2})ds\\
&\leq C\E\int_{0}^{T}(|Y(s)|^{p-2}+|D(\bar{Y}(s-\tau))|^{p-2})(1+|\bar Y(s)|^{2}+|\bar{Y}(s-\tau)|^{2})ds\\
&\leq C\E\int_{0}^{T}(1+|Y(s)|^{p}+|\bar{Y}(s)|^{p}+|\bar{Y}(s-\tau)|^{p})ds\leq C+C\int_{0}^{T}\E(\sup_{0\le u\le s}|Y(u)|^p)ds.
\end{split}
\end{equation}
Using Young's inequality, we derive that
\begin{equation}\label{H3}
\begin{split}
&\E(\sup_{0\leq t\leq T}\bar H_{2}(t))\\
&\leq  C\E\int_{0}^{T}|Y(s)-D(\bar{Y}(s-\tau))|^{p-2}\langle Y(s)-\bar{Y}(s),b_{h}(\bar{Y}(s),\bar{Y}(s-\tau))\rangle ds\\
&\leq C\E\big(\sup_{0\leq t\leq T}|Y(t)-D(\bar{Y}(t-\tau))|^{p-2}\times\big[(\int_{0}^{T}|Y(s)-\bar{Y}(s)||b_{h}(\bar{Y}(s),\bar{Y}(s-\tau))| ds)^{1/2}\big]^{2}\big)\\
&\leq \frac{1}{4}\E\sup_{0\leq t\leq T}|Y(t)-D(\bar{Y}(t-\tau))|^{p}+C\E\int_{0}^{T}|Y(s)-\bar{Y}(s)|^{p/2}|b_{h}(\bar{Y}(s),\bar{Y}(s-\tau))|^{p/2}ds\\
&\leq  \frac{1}{4}\E\sup_{0\leq t\leq T}|Y(t)-D(\bar{Y}(t-\tau))|^{p}+C,
\end{split}
\end{equation} 
where we have applied result from obtained Lemma \ref{Lemma1} to the second term above. 

By the BDG inequality again, we have
\begin{equation}\label{H4}
\begin{split}
&\E ( \sup_{0\leq t\leq T}\bar H_4(t))\\
&\E(\sup_{0\leq t\leq T}p\int_{0}^{t}|Y(s)-D(\bar{Y}(s-\tau))|^{p-2}\langle Y(s)-D(\bar{Y}(s-\tau)),\sigma(\bar{Y}(s),\bar{Y}(s-\tau))dB(s)\rangle)\\
&\leq C\E\bigg(\int_{0}^{T}\big(|Y(s)-D(\bar{Y}(s-\tau))|^{2p-2}\|\sigma(\bar{Y}(s),\bar{Y}(s-\tau))\|^{2}\big)ds\bigg)^{1/2}\\
&\leq C\E\bigg(\sup_{0\leq t\leq T}|Y(t)-D(\bar{Y}(t-\tau))|^{p-1}\big(\int_{0}^{T}\|\sigma(\bar{Y}(s),\bar{Y}(s-\tau))\|^{2}ds\big)^{1/2}\bigg)\\
&\leq \frac{1}{4}\E(\sup_{0\leq t \leq T}|Y(t)-D(\bar{Y}(t-\tau))|^{p})+C+C\int_{0}^{T}\E(\sup_{0\leq u\leq s}|Y(u)|^{p})ds.
\end{split}
\end{equation}
Substituting \eqref{H3} and \eqref{H4} into \eqref{EstimatedIntervalIto} and using  \eqref{Lemma301},  for  $p = 4,$ we can derive that
\begin{equation*}
\begin{split}
&\E(\sup_{0\leq t \leq T}|Y(t)|^{p}) \le C+ C\E(\sup_{0\leq t \leq T}|Y(t)-D(\bar{Y}(t-\tau))|^{p})\\
&\leq C+C\int_0^T\E(\sup_{0\le u\le s}|Y(u)|^p)ds.
\end{split}
\end{equation*}
The required assertion follows from an application of the Grownwall inequality. By repeating the same procedure, the desired result \eqref{UpperBuondLemma3} can be obtained by induction.
\hfill $\Box$

\subsection{Proof of Main Results}
In this subsection, we give proofs for the main results.

\textbf{Proof of Theorem \ref{ConvergenceRate}:} For any $0\leq t\leq T,$ by \cite[Lemma 6.4.1]{{Mao}} as well as the assumption {\bf(A2)}, we have, for $p\geq2$ and $\e>0,$
\begin{equation}\label{eq28}
\begin{split}
&|X(t)-Y(t)|^{p} = | X(t)-Y(t)-D(X(t-\tau))+D(\bar{Y}(t-\tau))\\
&~~~~~~~~~~~~~~~~~~~~~+D(X(t-\tau))-D(\bar{Y}(t-\tau))|^{p}\\
&\leq \bigg[1+\e^{\frac{1}{p-1}}\bigg]^{p-1}\bigg(\frac{|D(X(t-\tau))-D(\bar{Y}(t-\tau))|^{p}}{\e}\\
&~~~~~~~~~~~~~~~~~~~~~~+|X(t)-D(X(t-\tau))-Y(t)+D(\bar{Y}(t-\tau))|^{p}\bigg)\\
&\leq \bigg[1+\e^{\frac{1}{p-1}}\bigg]^{p-1}\bigg(\frac{\kappa^p|X(t-\tau)-\bar{Y}(t-\tau)|^{p}}{\e}\\
&~~~~~~~~~~~~~~~~~~~~~~+|X(t)-D(X(t-\tau))-Y(t)+D(\bar{Y}(t-\tau))|^{p}\bigg).
\end{split}
\end{equation}
Letting $\e=[\frac{\kappa}{1-\kappa}]^{p-1},$ \eqref{eq28} yields that
\begin{equation}
\begin{split}
|X(t)-Y(t)|^{p}&\le \kappa|X(t-\tau)-\bar{Y}(t-\tau)|^p\\
&+\frac{1}{(1-\kappa)^{p-1}}|X(t)-D(X(t-\tau))-Y(t)+D(\bar{Y}(t-\tau))|^{p}\\
&\leq \kappa|X(t)-\bar{Y}(t)|^pI_{[-\tau\leq t\leq 0]}(t)+\kappa|X(t)-\bar{Y}(t)|^pI_{[0\leq t\leq T]}(t)\\
&+\frac{1}{(1-\kappa)^{p-1}}|X(t)-D(X(t-\tau))-Y(t)+D(\bar{Y}(t-\tau))|^{p}.
\end{split}
\end{equation}
This, together with ({\bf A5}),  implies
 \begin{equation}\label{RateMainTheorem02}
 \begin{split}
 \E\big(\sup_{0\leq t\leq T}|X(t)-Y(t)|^{p}\big)&\leq \frac{1}{(1-\kappa)^{p}}\E\big(\sup_{0\leq t\leq T}|X(t)-D(X(t-\tau))\\
 &-Y(t)+D(\bar{Y}(t-\tau))|^{p}\big)+Ch^p.
 \end{split}
 \end{equation}
By It\^o's formula we can show that
\begin{equation}\label{GlobalProof01}
\begin{split}
&|X(t)-D(X(t-\tau))-Y(t)+D(\bar{Y}(t-\tau))|^{p}\\
&\leq p\int_{0}^{t}|X(s)-D(X(s-\tau))-Y(s)+D(\bar{Y}(s-\tau))|^{p-2}\\
&\times\langle X(s)-D(X(s-\tau))-Y(s)+D(\bar{Y}(s-\tau)),b(X(s),X(s-\tau))-b_{h}(\bar{Y}(s),\bar{Y}(s-\tau))\rangle ds\\
&+\frac{p(p-1)}{2}\int_{0}^{t}|X(s)-D(X(s-\tau))-Y(s)+D(\bar{Y}(s-\tau))|^{p-2}\\
&\times\|\sigma(X(s),X(s-\tau))-\sigma(\bar{Y}(s),\bar{Y}(s-\tau))\|^{2}ds\\
&+p\int_{0}^{t}|X(s)-D(X(s-\tau))-Y(s)+D(\bar{Y}(s-\tau))|^{p-2}\\
&\times\langle X(s)-D(X(s-\tau))-Y(s)+D(\bar{Y}(s-\tau)),\\
&~~~~~~~(\sigma(X(s),X(s-\tau))-\sigma(\bar{Y}(s),\bar{Y}(s-\tau)))dB(s)\rangle\\
&= p\int_{0}^{t}|X(s)-D(X(s-\tau))-Y(s)+D(\bar{Y}(s-\tau))|^{p-2}\\
&\K\langle X(s)-D(X(s-\tau))-Y(s)+D(\bar{Y}(s-\tau)),b(X(s),X(s-\tau))-b(Y(s),\bar{Y}(s-\tau))\rangle ds \\
&+p\int_{0}^{t}|X(s)-D(X(s-\tau))-Y(s)+D(\bar{Y}(s-\tau))|^{p-2}\\
&\K\langle X(s)-D(X(s-\tau))-Y(s)+D(\bar{Y}(s-\tau)),b(Y(s),\bar{Y}(s-\tau))-b(\bar{Y}(s),\bar{Y}(s-\tau))\rangle ds \\
& +p\int_{0}^{t}|X(s)-D(X(s-\tau))-Y(s)+D(\bar{Y}(s-\tau))|^{p-2}\\
&\K\langle X(s)-D(X(s-\tau))-Y(s)+D(\bar{Y}(s-\tau)),b(\bar{Y}(s),\bar{Y}(s-\tau))-b_{h}(\bar{Y}(s),\bar{Y}(s-\tau))\rangle ds\\
&+\frac{p(p-1)}{2}\int_{0}^{t}|X(s)-D(X(s-\tau))-Y(s)+D(\bar{Y}(s-\tau))|^{p-2}\\
&\times||\sigma(X(s),X(s-\tau))-\sigma(\bar{Y}(s),\bar{Y}(s-\tau))||^{2}ds\\
&+p\int_{0}^{t}|X(s)-D(X(s-\tau))-Y(s)+D(\bar{Y}(s-\tau))|^{p-2}\\
&\times\langle X(s)-D(X(s-\tau))-Y(s)+D(\bar{Y}(s-\tau)),\\
&~~~~~~~~(\sigma(X(s),X(s-\tau))-\sigma(\bar{Y}(s),\bar{Y}(s-\tau)))dB(s)\rangle\\
& =: I_{1}(t)+I_{2}(t)+I_{3}(t)+I_{4}(t)+I_{5}(t).
\end{split}
\end{equation}
By assumption {\bf(A2)} and  {\bf(A4)},  we have
\begin{equation}\label{I1}
\begin{split}
&\E(\sup_{0\leq t\leq T}(I_{1}(t)))\leq C\E\int_{0}^{T}|X(s)-D(X(s-\tau))-Y(s)+D(\bar{Y}(s-\tau))|^{p-2}\\
&\K (|X(s)-Y(s)|^{2}+|X(s-\tau)-\bar{Y}(s-\tau)|^{2})ds\\
&\leq C\E\int_{0}^{T}\big(|X(s)-Y(s)|^{p-2}+|D(X(s-\tau))-D(\bar{Y}(s-\tau))|^{p-2}\big)\\
&\K (|X(s)-Y(s)|^{2}+|X(s-\tau)-\bar{Y}(s-\tau)|^{2})ds\\
&\leq C\int_{0}^{T}\E(\sup_{0\leq u\leq s}|X(u)-Y(u)|^{p})ds+C\int_{0}^{T}\E(\sup_{0\leq u\leq s}|X(u)-Y(u)+Y(u)-\bar{Y}(u)|^{p})ds\\
&+C\int_{-\tau}^{0}\E(\sup_{-\tau\leq \theta\leq 0}|X(\theta)-\bar{Y}(\theta)|^{p})ds\leq C\int_{0}^{T}\E(\sup_{0\leq u\leq s}|X(u)-Y(u)|^{p})ds+Ch^{p}+Ch^{p/2}.
\end{split}
\end{equation}
In the same way as in \eqref{I1}, we can estimate  $I_{4}(t)$  such that 
\begin{equation}\label{I4}
\E(\sup_{0\leq t\leq T}(I_{4}(t)))\leq C\int_{0}^{T}\E(\sup_{0\leq u\leq s}|X(u)-Y(u)|^{p})ds+Ch^{p}+Ch^{p/2}.
\end{equation}
Using assumption {\bf(A3)}, H\"older's inequality and  \eqref{EstimateIntervalError}, we arrive at
\begin{equation}\label{I2}
\begin{split}
\E(\sup_{0\leq t\leq T}I_{2}(t))&\leq \E\bigg(\int_{0}^{T}|X(s)-D(X(s-\tau))-Y(s)+D(\bar{Y}(s-\tau))|^{p-1}\\
&\K L(1+|\bar{Y}(s)|^{l}+|Y(s)|^{l}+2|\bar{Y}(s-\tau)|^{l})(|Y(s)-\bar{Y}(s)|)ds\bigg)\\
&\leq\E\bigg(\sup_{0\leq s\leq T}|X(s)-D(X(s-\tau))-Y(s)+D(\bar{Y}(s-\tau))|^{p-1}\\
&\K \int_{0}^{T}L(1+|\bar{Y}(s)|^{l}+|Y(s)|^{l}+2|\bar{Y}(s-\tau)|^{l})(|Y(s)-\bar{Y}(s)|)ds\bigg)\\
&\leq \frac{1}{4}\E(\sup_{0\leq s\leq T}|X(s)-D(X(s-\tau))-Y(s)+D(\bar{Y}(s-\tau))|^{p})\\
&+C\E\bigg[\int_{0}^{T}(1+|\bar{Y}(s)|^{l}+|Y(s)|^{l}+|\bar{Y}(s-\tau)|^{l})(|Y(s)-\bar{Y}(s)|)ds\bigg]^{p}\\
&\leq  \frac{1}{4}\E(\sup_{0\leq s\leq T}|X(s)-D(X(s-\tau))-Y(s)+D(\bar{Y}(s-\tau))|^{p})\\
&+ C\bigg[\int_{0}^{T}\E(1+|\bar{Y}(s)|^{l}+|Y(s)|^{l}+|\bar{Y}(s-\tau)|^{l})^{p}(|Y(s)-\bar{Y}(s)|)^{p}ds\bigg]\\
&\leq \frac{1}{4}\E(\sup_{0\leq s\leq T}|X(s)-D(X(s-\tau))-Y(s)+D(\bar{Y}(s-\tau))|^{p})\\
&+C[\int_{0}^{T}(\E(1+|\bar{Y}(s)|^{l}+|Y(s)|^{l}+|\bar{Y}(s-\tau)|^{l})^{2p})^{1/2}(\E|Y(s)-\bar{Y}(s)|^{2p})^{1/2}ds]\\
&\leq \frac{1}{4}\E(\sup_{0\leq s\leq T}|X(s)-D(X(s-\tau))-Y(s)+D(\bar{Y}(s-\tau))|^{p})+Ch^{p/2}.
\end{split}
\end{equation}

We now estimate $I_{3},$ 
\begin{equation}\label{I3}
\begin{split}
\E(\sup_{0\leq t\leq T}I_{3}(t))& \leq p\E\int_{0}^{T}|X(s)-D(X(s-\tau))-Y(s)+D(\bar{Y}(s-\tau))|^{p-1}\\
&\K|b(\bar{Y}(s),\bar{Y}(s-\tau))-b_{h}(\bar{Y}(s),\bar{Y}(s-\tau))|ds\\
&\leq\frac{1}{4}\E(\sup_{0\leq s\leq T}|X(s)-D(X(s-\tau))-Y(s)+D(\bar{Y}(s-\tau))|^{p})\\
&+C\E\int_{0}^{T}|b(\bar{Y}(s),\bar{Y}(s-\tau))-b_{h}(\bar{Y}(s),\bar{Y}(s-\tau))|^{p}ds\\
&\leq \frac{1}{4}\E(\sup_{0\leq s\leq T}|X(s)-D(X(s-\tau))-Y(s)+D(\bar{Y}(s-\tau))|^{p})+Ch^{\alpha p},
\end{split}
\end{equation}
where we have used the following fact
\begin{equation}\label{RateProof03}
\begin{split}
& \E\int_{0}^{T}|b(\bar{Y}(s),\bar{Y}(s-\tau))-b_{h}(\bar{Y}(s),\bar{Y}(s-\tau))|^{p}ds\\
&\leq  h^{\alpha p} \E\bigg[\int_{0}^{T}\frac{|b(\bar{Y}(s),\bar{Y}(s-\tau))|^{2p}}{(1+h^{\alpha}|b(\bar{Y}(s),\bar{Y}(s-\tau))|)^{p}}ds\bigg]\\
&\leq h^{\alpha p} \int_{0}^{T}\E\bigg[|L(1+|\bar{Y}(s)|^{l}+|\bar{Y}(s-\tau)|^{l})(|\bar{Y}(s)|+|\bar{Y}(s-\tau)|)+C|^{2p}\bigg]ds\\
&\leq Ch^{\alpha p}.
\end{split}
\end{equation}
Moreover, the BGD inequality yields that
\begin{equation}\label{I5}
\begin{split}
&\E(\sup_{0\leq t\leq T}I_{5}(t))\leq C\E\big(\int_{0}^{T}|X(s)-D(X(s-\tau))-Y(s)+D(\bar{Y}(s-\tau))|^{2p-2}\\
&\K\|\sigma(X(s),X(s-\tau))-\sigma(\bar{Y}(s),\bar{Y}(s-\tau))\|^{2}ds\big)^{1/2}\\
&\leq \frac{1}{4}\E(\sup_{0\leq s\leq T}|X(s)-D(X(s-\tau))-Y(s)+D(\bar{Y}(s-\tau))|^{p})\\
&+C\int_{0}^{T}\E(|X(s)-\bar Y(s)|^{p}+|X(s-\tau)-\bar{Y}(s-\tau)|^{p})ds\\
&\leq \frac{1}{4}\E(\sup_{0\leq s\leq T}|X(s)-D(X(s-\tau))-Y(s)+D(\bar{Y}(s-\tau))|^{p})\\
&+C\int_{0}^{T}\E(\sup_{0\leq u\leq s}|X(u)-Y(u)|^{p})ds+Ch^{p}.
\end{split}
\end{equation}
Substituting (\ref{I1}), (\ref{I4}), (\ref{I2}), (\ref{I3}) and (\ref{I5}) into (\ref{RateMainTheorem02}).
\begin{equation}\label{RateMainTheorem03}
 \begin{split}
 \E\big(\sup_{0\leq t\leq T}|X(t)-Y(t)|^{p}\big)&\leq Ch^{\alpha p}+Ch^{p/2}+Ch^p+C\int_{0}^{T}\E(\sup_{0\leq u\leq s}|X(u)-Y(u)|^{p})ds\\
 &\leq Ch^{\a p}+C\int_{0}^{T}E(\sup_{0\leq u\leq s}|X(u)-Y(u)|^{p})ds.
 \end{split}
 \end{equation}
The desired result follows by the Grownwall inequality and  the fact $\alpha\in(0,1/2]$.  \hfill
 $\Box$

In the sequel, we give a proof for Theorem \ref{Convergence} which the monotonicity condition {\bf (A4)} is replaced by the local one {\bf (A3)}. The techniques of the proof have been developed in Higham, Mao and
Stuart \cite{HMS} where they showed the strong convergence of the EM method
for the SDE  under a local Lipschitz condition. We
therefore only give the outline of the proof here.

\textbf{Proof of Theorem \ref{Convergence}:}  For every $R >0,$ we define the  stopping times
\begin{align}\label{StopppingTimes}
	\tau_{R}:= \inf\{t\geq 0: |X(t)| \geq R\}, \quad \rho_{R}:= \inf\{t\geq 0: |Y(t)| \geq R\},
\end{align}
and  denote that $\theta_{R} = \tau_{R}\wedge  \rho_{R},$ and
$
e(t)=X(t)-Y(t).
$
By the virtue of Young's inequality, for $q>p$ and $\eta >0$ we have
\begin{equation*}
\begin{split}
&\E\bigg[\sup_{0\leq t\leq T}|e(t)|^{p}\bigg] \leq \E\bigg[\sup_{0\leq t\leq T}|e(t)|^{p}I_{\{\text{$\tau_{R}\leq T$ or $\rho_{R}$} \leq T\}}\bigg]+\E\bigg[\sup_{0\leq t\leq T}|e(t\wedge \theta_{R})|^{p}\bigg]\\
&=\E\bigg[\bigg(\eta^{p/q}\sup_{0\leq t\leq T}|e(t)|^{p}\bigg)\bigg(\eta^{-p/q}I_{\{\text{$\tau_{R}\leq T$ or $\rho_{R}$} \leq T\}}\bigg)\bigg] +\E\bigg[\sup_{0\leq t\leq T}|e(t\wedge \theta_{R})|^{p}\bigg]\\
&\leq \frac{p\eta}{q}\E\bigg[\sup_{0\leq t \leq T}|X(t)-Y(t)|^{q}\bigg]+\frac{q-p}{q\eta^{p/(q-p)}}P(\text{$\tau_{R}\leq T$ or $\rho_{R} \leq T$}) \\
&+\E\bigg[\sup_{0\leq t\leq T}|e(t\wedge \theta_{R})|^{p}\bigg]\\
&\leq \frac{p\eta}{q}2^{q}\E\bigg\{\bigg[\sup_{0\leq t \leq T}|X(t)|^{q}\bigg]+\bigg[\sup_{0\leq t \leq T}|X(t)|^{q}\bigg]\bigg\}+\frac{q-p}{q\eta^{p/(q-p)}}\bigg\{\E\bigg[\frac{|X(\tau_{R})|}{R^{p}}\bigg]+\E\bigg[\frac{|Y(\rho_{hR})|}{R^{p}}\bigg]\bigg\}\\
&+\E\bigg[\sup_{0\leq t\leq T}|e(t\wedge \theta_{R})|^{p}\bigg]\\
&\leq \frac{p\eta}{q}2^{q}C+\frac{q-p}{q\eta^{p/(q-p)}R^p}C+\E\bigg[\sup_{0\leq t\leq T}|e(t\wedge \theta_{R})|^{p}\bigg].\\
\end{split}
\end{equation*}
In the similar way as Theorem \ref{ConvergenceRate} was proved, we can show that
\begin{equation}
\E\bigg[\sup_{0\leq t\leq T}|e(t\wedge \theta_{R})|^{p}\bigg]\le C_R h+o(h).
\end{equation}
Finally, given an $\epsilon >0,$ there exist some $\eta$ small enough that
\begin{equation*}
\frac{p\eta}{q}2^{q}C < \frac{\epsilon}{3},
\end{equation*}
choose $R$ large enough that
\begin{equation*}
\frac{q-p}{q\eta^{p/(q-p)}R^p}2C < \frac{\epsilon}{3},
\end{equation*}
and $h$ small enough that
\begin{equation*}
\E[\sup_{t_0\leq t\leq T}|e(t\wedge \theta_R)|^{p}] < \frac{\epsilon}{3},
\end{equation*}
Hence we obtain
\begin{equation*}
\E\bigg[\sup_{0\leq t\leq T}|X(t)-Y(t)|^p\bigg]< \epsilon.
\end{equation*}
\hfill $\Box$

\section{Tamed EM Method of NSDDEs driven by Pure Jump Processes}
In this section, we investigate NSDDEs driven by pure jump processes. Similar to NSDDEs driven by Brownian motion, most  NSDDEs driven by pure jumps have no explicit solutions. Therefore, it is important to investigate the numerical approximation of NSDDEs driven by pure jumps.

We need to introduce more notation. Let $\mathcal{B}(\R)$ be the Borel $\sigma$-algebra on $\R$ and $\lambda(dx)$ be a $\sigma$-finite measure defined on $\mathcal{B}(\R).$ Let $D(\varrho)\subset(0,\infty)$ be a countable set, and let $\varrho=(\varrho(t), t\in D(\mathbb{\varrho}))$ be a stationary $\mathcal{F}_{t}$-Poisson point process on $\R$ with characteristic measure $\lambda(\cdot).$ Denote that by $N(dt,du)$ the Poisson counting measure associated with $\varrho,$ i.e., $N(t,U)= \sum_{s\in D(\varrho),s\leq t}I_{U}(\varrho(s))$ for $U\in \mathcal{B}(\R\setminus\{0\}).$ Let $\tilde{N}(dt,du) = N(dt,du)-dt\lambda(du)$ be the compensated Poisson measure associated with $N(dt,du).$ In what follows, we further assume that $\int_{U}|u|^p\lambda(du) <\infty$ for any $p\geq 1.$ 
In this section, we consider the following NSDDE driven by pure jump processes on $\R^{n}:$
\begin{equation}\label{NSDDEJump}
d[x(t)-G(x(t-\tau))]  = f(x(t),x(t-\tau))dt+\int_{U}g(x(t-),x((t-\tau)-),u)\tilde{N}(du,dt),~t\ge0,
\end{equation}
with initial data $\{x(\theta):-\tau\leq \theta\leq 0\} = \xi \in L^p_{\mathcal{F}_{0}}([-\tau,0];\mathbb{R}^{n}), p\ge 2,$ $x(t-):= \lim_{s\uparrow t}x(t),$
where $G:\R^n\mapsto\R^n$, $f:\R^n\times\R^n\mapsto\R^n$, and
$g:\R^n\times\R^n\times U\mapsto\R^n$ are  Borel measurable. 

Again, we assume that the step size $h\in (0, 1)$ be fraction of two positive rational numbers $\tau$ and $T,$  so that there exist two positive integers $M, \bar{M}$ such that $h=T/M=\tau/\bar{M}.$

For the future use, we assume: 
\begin{description}
\item[(B1)]There exists a positive constant $K_{1}$ such that
\begin{equation}\label{WeakCorecivityJump}
2\langle x-G(y),f(x,y)\rangle \leq K_{1}(1+|x|^{2}+|y|^2)  \mbox{ and }\int_{U}|g(x,y,u)|^{p}\lambda(du)\leq K_{1}(1+|x|^{p}+|y|^p)
\end{equation}
for $\forall$ $x,y \in \mathbb{R}^n,  p\ge 2.$ 
\item[(B2)] $G(0)=0$ and there exists a constant $\kappa \in (0,1)$ such that
\begin{eqnarray}\label{NeutralBoundJump}
|G(x)-G(\bar{x})| \leq \kappa|x-\bar{x}| \mbox{ for all } x,y \in \mathbb{R}^n.
\end{eqnarray}
\item[(B3)] For any $R>0,$ there exist two positive constants $\tilde{K}_{R}$ and $K_{R}$ such that for all  $|x|\vee|y|\vee|\bar{x}|\vee|\bar{y}| \leq R,   p\ge 2$

\begin{equation}\label{LocalLipschitzIIJump}
\begin{split}
&2\langle x-G(y)-\bar{x}+G(\bar{y}), f(x,y)-f(\bar{x},\bar{y})\rangle
\leq \tilde{K}_{R} (|x-\bar{x}|^{2}+|y-\bar{y}|^2), \\
&\int_{U}|g(x,y,u)-g(\bar{x},\bar{y},u)|^{p}\lambda(du)\le \tilde{K}_{R} (|x-\bar{x}|^{p}+|y-\bar{y}|^p),
\end{split}
\end{equation}
and 
\begin{equation*}
\sup_{|x|\leq R, |y|\leq R}|f(x,y)| \leq K_{R}.
\end{equation*}
\item[(B4)] There exist two positive constants $l$ and $L$ such that for all $x, y, \bar x, \bar y \in \R^n, p\ge 2$
\begin{equation*}
\begin{split}
&2\langle x-G(y)-\bar{x}+G(\bar{y}), f(x,y)-f(\bar{x},\bar{y})\rangle
\leq L(|x-\bar{x}|^2 + |y-\bar{y}|^2)\\
&\int_{U}|g(x,y,u)-g(\bar{x},\bar{y},u)|^{p}\lambda(du)\le  L(|x-\bar{x}|^p + |y-\bar{y}|^p)
\end{split}
\end{equation*}
and
\begin{equation*}
|f(x,y)-f(\bar{x},\bar{y})|\leq L(1+|x|^{l}+|y|^{l}+|\bar{x}|^{l}+|\bar{y}|^{l})(|x-\bar{x}|+|y-\bar{y}|).
\end{equation*}

\item[(B5)] For every $p >0,$ there exists a positive integer $K,$ such that
$$\E\|\xi(t)-\xi(s)\|^{p}\leq K|t-s|^p, \mbox{ for any } s,t\in[-\tau,0].$$ 
\end{description}
\begin{remark}\label{solJump}
{\rm Under assumptions {\bf(B1)}, {\bf(B2)} and {\bf(B4)} , the NSDDE \eqref{NSDDEJump} with the initial data $x(0)=\xi $ satisfying $\E\|\xi\|^2<\infty$ has a pathwise unique strong solution. The details of proof  can be found in  \cite[Theorem 170]{Situ}. If the  condition {\bf(B4)} is replaced by the condition {\bf(B3)}, the existence and uniqueness theorem still holds.}
\end{remark}
Set
$$
f_h(x, y)=\frac{f(x, y)}{1+h^\alpha|f(x, y)|}, \, \, \alpha\in(0,1/2).
$$
Similarly to the Brownian motion case, the discrete-time tamed EM scheme associated with (\ref{NSDDEJump}) can be defined as following:
For every integer $n = -\bar{M},\cdots,0,$ $z_{h}^{(n)} = \xi(nh).$
For every integer $n = 1, \cdots, M-1,$
\begin{equation}\label{DEulerschemeStepJump}
z_{h}^{(n+1)}-G(z_{h}^{(n+1-\bar{M})})= z_{h}^{(n)}-G(z_{h}^{(n-\bar{M})})+f_{h}(z_{h}^{(n)},z_{h}^{(n-\bar{M})})h+g(z_{h}^{(n)},z_{h}^{(n-\bar{M})},u)\Delta \tilde{N}^{n}_{h},
\end{equation}
where $\Delta \tilde{N}^{n}_{h}:= \tilde{N}((n+1)h,U)-\tilde{N}(nh,U),$ the increment of compensated Poisson process.
We may rewrite  the discrete tamed EM scheme  as follows:
\begin{equation}\label{DEulerschemeSumJump}
\begin{split}
z_{h}^{(n+1)}=&G(z_{h}^{(n+1-\bar{M})})+ \xi(0)-G(\xi(-\tau))+\sum_{i=1}^{n}f_{h}(z_{h}^{(i)},z_{h}^{(i-\bar{M})})h\\
&+\sum_{i=1}^{n}g(z_{h}^{(i)},z_{h}^{(i-\bar{M})},u)\Delta \tilde{N}^{i}_{h}.
\end{split}
\end{equation}
For $t\in[nh,(n+1)h),$ denote that $\bar{z}(t): = z_{h}^{(n)},$ and then $\bar{z}(t-\tau) = z_{h}^{(n-\bar{M})}.$ 
It is more convenient to define the continuous-time tamed EM approximate solution $z(t)$ associated with (\ref{NSDDEJump}) as below:

For any $\theta\in[-\tau,0],$ $z(\theta) = \xi(\theta).$ For any $t \in [0,T],$
\begin{equation}\label{CEulerschemeJump}
\begin{split}
z(t)=& G(\bar{z}(t-\tau))+\xi(0)-G(\xi(-\tau))+\int_{0}^{t}f_{h}(\bar{z}(s),\bar{z}(s-\tau))ds\\
&+\int_{0}^{t}\int_{U}g(\bar{z}(s-),\bar{z}((s-\tau)-),u)\tilde{N}(du,ds).
\end{split}
\end{equation}
Since for any $t>0$ there exists a positive integer $n,$ $0\leq n\leq M-1,$ such that $t \in[nh,(n+1)h),$ we have
\begin{equation}\label{CEulerschemeModificationJump}
\begin{split}
z(t) 
= z(nh)+\int_{nh}^{t}f_{h}(\bar{z}(s),\bar{z}(s-\tau))ds+\int_{nh}^{t}\int_{U}g(\bar{z}(s-),\bar{z}((s-\tau)-),u)\tilde{N}(du,ds).
\end{split}
\end{equation}
Clearly, continuous-time tamed EM approximate solution $z(t)$ coincides with the discrete-time tamed approximation solution $\bar{z}(t)$ at the grid points $t = nh,$ i.e. $\bar{z}(t) = z^{n}_{h}=z(nh).$ 

We now state our main results of this section.
\begin{theorem}\label{ConvergenceRateJump}
	Assume that {\bf(B1)}, {\bf(B2)}, {\bf(B4)} and {\bf(B5)} hold. Assume also that  $p \geq 2$ and $\alpha \in (0, 1/p)$, then  the tamed EM scheme $(\ref{CEulerschemeJump})$ converges to the exact solution of $(\ref{NSDDEJump})$ such that
	\begin{equation}\label{MainResultConvergenceRateJump}
	\E\Big[\sup_{0 \leq t \leq T}|z(t)-x(t)|^{p}\Big] \leq Ch^{\gamma},
	\end{equation}
	where $\g = 1/2 \wedge \alpha p.$
\end{theorem}
and
\begin{theorem}\label{ConvergenceJump}
	Assume that that {\bf(B1)}-{\bf(B3)} and {\bf(B5)} hold. Assume also that  $p \geq 2$ and $\alpha \in (0, 1/p)$,  then the tamed EM scheme $(\ref{CEulerschemeJump})$ converges to the exact solution of $(\ref{NSDDEJump})$ such that
	\begin{equation}\label{MainResultConvergenceJump}
	\lim_{h \to 0}\E\Big[\sup_{0 \leq t \leq T}|z(t)-x(t)|^{p}\Big] = 0. 
	\end{equation}
\end{theorem}

\subsection{Boundedness of Moments}

In order to show our main results, we need the  following inequality \cite[Theorem 1]{{MAR}}.
\begin{lemma}\label{BJ}
Let $\varphi: \R_{+}\times \to \R^{n}$ and assume that 
\begin{equation*}
\int_{0}^{t}\int_{U}\E|\varphi(s,u)|^{2}\lambda(du)ds<\infty, \quad t\geq 0 \quad
\mbox{a.s.}.
\end{equation*}
Then there exist a constant $C_{p} >0$ such that for any $p>0$
\begin{equation*}
\begin{split}
\E\big(\sup_{0\leq s\leq t}|\int_{0}^{s}\int_{U}\varphi(r-,u)\tilde{N}(dr,du)|^{p}\big)&\leq C_{p}\bigg[\E\bigg(\int_{0}^{t}\int_{U}|\varphi(s,u)|^{2}\lambda(du)d s\bigg)^{p/2}\\
&+\E\int_{0}^{t}\int_{U}|\varphi(s,u)|^{p}\lambda(du)ds\bigg].
\end{split}
\end{equation*}	
\end{lemma}
It is known that if $1\leq p\leq 2,$ then the second term on the right hand side can be eliminated.

\begin{lemma}\label{yl3Jump}
Consider the continuous-time tamed EM scheme given by equation $(\ref{CEulerschemeModificationJump}).$  Assume that $p \geq 2, \alpha\in(0,1/p),$  and 
\begin{equation}\label{AssumptionLemma1Jump}
\sup_{0\leq t\leq T}\E(|z(t)|^{p}) \le C.
\end{equation}
Assume also that {\bf (B1)} holds, then it holds that 
\begin{equation}\label{EstimateIntervalErrorJump}
\E\bigg[\sup_{0\leq n\leq M-1}\sup_{nh\leq t \leq (n+1)h}|z(t)-z(nh)|^p\bigg]\leq Ch,
\end{equation}
and
\begin{equation}\label{ErrorWithDriftJump}
\E\bigg[\sup_{0\leq n\leq M-1}\sup_{nh\leq t \leq  (n+1)h}|z(t)-z(nh)|^p|f_{h}(\bar{z}(t),\bar{z}(t-\tau))|^{p}\bigg] \leq C.
\end{equation}
\end{lemma}
\textbf{Proof}:
By the definition of $z(t)$, we can write that
\begin{equation*}
\begin{split}
\E\bigg[\sup_{nh\leq t \leq(n+1)h}|z(t)-z(nh)|^p\bigg] &= \E\bigg[\sup_{nh\leq t \leq(n+1)h}\bigg|\int_{nh}^{t}f_{h}(\bar{z}(s),\bar{z}(s-\tau))ds\\
&+\int_{nh}^{t}\int_{U}g(\bar{z}(s-),\bar{z}((s-\tau)-),u)\tilde{N}(du,ds)\bigg|^p\bigg].
\end{split}
\end{equation*}
Therefore, due to H\"older's inequality,
\begin{equation}\label{Lemma1ObservationJump}
\begin{split}
&\E\bigg[\sup_{nh\leq t \leq (n+1)h}|z(t)-z(nh)|^p\bigg]\\
&\leq 2^{p-1} h^{p-1}\E\bigg[\int_{nh}^{(n+1)h}|f_{h}(\bar{z}(s),\bar{z}(s-\tau))|^{p}ds\bigg]\\
&+2^{p-1}\E\bigg[\sup_{nh\leq t \leq (n+1)h}\bigg|\int_{nh}^{t}\int_{U}g(\bar{z}(s-),\bar{z}((s-\tau)-),u)\tilde{N}(du,ds)\bigg|^{p}\bigg].
\end{split}
\end{equation}
Using the Lemma \ref{BJ} and  $(\ref{AssumptionLemma1Jump}),$ for some $p \geq 2$ we have
\begin{equation*}
\begin{split}
&\E\bigg[\sup_{nh\leq t \leq (n+1)h}\left|\int_{nh}^{t}\int_{U}g(\bar{z}(s-),\bar{z}((s-\tau)-),u)\tilde{N}(du,ds)\right|^{p}\bigg]\\
&\leq C\E\int_{nh}^{(n+1)h}\int_{U}|g(\bar{z}(s),\bar{z}(s-\tau),u)|^{p}\lambda(du)ds\\
&\leq C\E\int_{nh}^{(n+1)h}(1+|\bar{z}(s)|^{p}+|\bar{z}(s-\tau)|^{p})ds\\
&\leq Ch.
\end{split}
\end{equation*}
This, together with  $|f_{h}(\bar{z}(s),\bar{z}(s-\tau),s)|\leq h^{-\alpha},$ yields
\begin{equation}\label{ayJump}
\E\bigg[\sup_{nh\leq t \leq (n+1)h}|z(t)-z(nh)|^p \bigg]\leq  2^{p-1}h^{(1-\alpha)p}+Ch \leq Ch.
\end{equation}
Hence $(\ref{EstimateIntervalErrorJump})$ holds.
Moreover,
\begin{equation}\label{ay1Jump}
\begin{split}
&\E\bigg[\sup_{nh\leq t \leq (n+1)h}|z(t)-z(nh)|^p|f_{h}(\bar{z}(s),\bar{z}(s-\tau))|^{p}\bigg] \\
&\leq \E\bigg[\sup_{nh\leq t \leq (n+1)h}|z(t)-z(nh)|^p\bigg]h^{-\alpha p} 
\leq Ch^{1-\alpha p}\leq C,
\end{split}
\end{equation}
as required. The proof is therefore complete. \hfill $\Box$

The  results of boundedness are given below:
\begin{lemma}\label{Boundedness}
Assume that {\bf(B1)} and {\bf(B2)} hold. Assume also   $p \geq 2, \alpha\in (0, 1/p)$,  then there exists a constant $C$ such that
	\begin{equation}\label{UpperBuondLemma3Jump}
	\E[\sup_{0\leq t \leq T}|x(t)|^{p}]\vee \E[\sup_{0\leq t \leq T}|z(t)|^{p}] \leq C.
	\end{equation}	
\end{lemma}
\textbf{Proof}:\quad
An application of the It\^o formula yields,
 \begin{equation*}
 \begin{split}
 &|x(t)-G(x(t-\tau))|^{p} = |\xi(0)-G(\xi(-\tau))|^{p}\\
 &+p\int_{0}^{t}|x(s)-G(x(s-\tau))|^{p-2}\langle x(s)-G(x(s-\tau)),f(x(s),x(s-\tau))\rangle ds\\
 &+\int_{0}^{t}\int_{U}|x(s)-G(x(s-\tau))+g(x(s),x(s-\tau),u)|^{p}-|x(s)-G(x((s-\tau))|^{p}\\
 &-p|x(s)-G(x(s-\tau))|^{p-2}\langle x(s)-G(x(s-\tau)),g(x(s),x(s-\tau),u)\rangle \lambda(du)ds\\
 &+\int_{0}^{t}\int_{U}p|x(s)-G(x(s-\tau))|^{p-2}\\
 &\K\langle x(s)-G(x(s-\tau)),g(x(s-),x((s-\tau)-),u)\rangle\tilde{N}(du,ds)\\
 &=: J_{1}(t)+J_{2}(t)+J_{3}(t)+J_{4}(t).
 \end{split}
 \end{equation*}
 Due to assumption {\bf(B1)}, we derive that
 \begin{equation}\label{J2}
 \begin{split}
 &\E(\sup_{0\leq t\leq T}J_{2}(t))\\
 &\leq C\E\int_{0}^{T}|x(s)-G(x(s-\tau))|^{p-2}(1+|x(s)|^{2}+|x(s-\tau)|^{2})ds\\
 &\leq C+C\int_{0}^{T}\E\big(\sup_{0\leq u\leq s}|x(u)|^{p}\big)ds.
 \end{split}
 \end{equation}
By the taylor expansion, we have
 \begin{equation}\label{Taylor}
 \begin{split}
  |x+y|^{p} - |x|^{p}-|x|^{p-2}|xy|
  \leq C(|x|^{p-2}|y|^{2}+|y|^{p}).
 \end{split}
 \end{equation} 
 This, together with {\bf(B1)} and  {\bf(B2)},  implies
 \begin{equation}\label{J3}
 \begin{split}
 &\E(\sup_{0\leq t\leq T}J_{3}(t))= \E\int_{0}^{T}\int_{U}|x(s)-G(x(s-\tau))+g(x(s),x(s-\tau),u)|^{p}\\
 &-|x(s)-G(x(s-\tau))|^{p}-p|x(s)-G(x(s-\tau))|^{p-2}\\
 &\K\langle x(s)-G(x(s-\tau)),g(x(s),x(s-\tau),u)\rangle \lambda(du)ds\\
 &\leq \E\int_{0}^{T}\int_{U}\big[|x(s)-G(x(s-\tau))|^{p-2}|g(x(s),x(s-\tau),u)|^{2}\\
 &+|g(x(s),x(s-\tau),u)|^{p}\big]\lambda(du)ds\\
 &\leq \E\int_{0}^{T}|x(s)-G(x(s-\tau))|^{p-2}(1+|x(s)|^{2}+|x(s-\tau)|^{2})\\
 &+(1+|x(s)|^{p}+|x(s-\tau)|^{p})ds\\
 &\leq C+ C\int_{0}^{T}\E\big(\sup_{0\leq u\leq s}|x(u)|^{p}\big)ds.
 \end{split}
 \end{equation}
   Using  Lemma \ref{BJ} and Young's inequality, one has
 \begin{equation}\label{J4}
 \begin{split}
 &\E(\sup_{0\leq t\leq T}J_{4}(t))\\
 &\leq C\E\bigg(\int_{0}^{T}\int_{U}|x(s)-G(x(s-\tau))|^{2p-2}|g(x(s),x(s-\tau),u)|^{2} \lambda(du)ds\bigg)^{1/2}\\
 &\leq C\E[(\sup_{0\leq t\leq T}|x(t)-G(x(s-\tau))|^{p-1})\K\big(\int_{0}^{T}\int_{U}|g(x(s),x(s-\tau),u)|^{2}\big)^{1/2}]\\
 &\leq \frac{1}{4}\E(\sup_{0\leq t\leq T}|x(t)-G(x(s-\tau))|^{p})+C\E\big(\int_{0}^{T}(1+|x(s)|^{2}+|x(s-\tau)|^{2})\big)^{p/2}\\
 &\leq \frac{1}{4}\E(\sup_{0\leq t\leq T}|x(t)-G(x(s-\tau))|^{p})+C+C\int_{0}^{T}\E\big(\sup_{0\leq u\leq s}|x(u)|^{p}\big)ds.
 \end{split}
 \end{equation}
Applying \cite[Lemma 6.4.4]{Mao} (also see \eqref{Lemma301}), we  derive that
 \begin{equation*}
 \E\big(\sup_{0\leq t\leq T}|x(t)|^{p}\big)\leq C+C\int_{0}^{T}\E\big(\sup_{0\leq u\leq s}|x(u)|^{p}\big)ds.
 \end{equation*}
 By the Gronwall inequality, we have
 \begin{equation*}
 \E\big(\sup_{0\leq t\leq T}|x(t)|^{p}\big)\leq C.
 \end{equation*}
The proof of boundedness of EM approximation is analogous to its Brownian motion counterpart, we first claim that there exists a constant $C>0$ such that:
\begin{equation}\label{H404}
\sup_{0\leq t\leq T}\E(|z(t)|^{2})\leq C.
\end{equation}
Using the It\^o formula, we have
\begin{equation}\label{H400}
\begin{split}
&|z(t)-G(\bar{z}(t-\tau))|^{2} = |\xi(0)-G(\xi(-\tau))|^{2}\\
&+2\int_{0}^{t}\langle \bar z(s)-G(\bar{z}(s-\tau)),f_{h}(\bar{z}(s),\bar{z}(s-\tau))\rangle ds\\
&+2\int_{0}^{t}\langle  z(s)-\bar{z},f_{h}(\bar{z}(s),\bar{z}(s-\tau))\rangle ds\\
&+\bigg(\int_{0}^{t}\int_{U}|z(s)-G(\bar{z}(s-\tau))+g(\bar{z}(s),\bar{z}(s-\tau),u)|^{2}-|z(s)-G(\bar{z}(s-\tau))|^{2}\\
&-2\langle z(s)-G(\bar{z}(s-\tau)),g(\bar{z}(s),\bar{z}(s-\tau),u)\rangle\bigg) \lambda(du)ds\\
&+\int_{0}^{t}\int_{U}2\langle z(s-)-G(\bar{z}((s-\tau)-)), g(\bar{z}(s-),\bar{z}((s-\tau)-),u)\rangle\tilde{N}(du,ds)\\
&=:  |\xi(0)-G(\xi(-\tau))|^{2}+ \bar{J}_{1}(t)+\bar{J}_{2}(t)+\bar{J}_{3}(t)+\bar{J}_{4}(t).
\end{split}
\end{equation}
By {\bf(B1)}, we compute
\begin{equation}\label{H401}
\begin{split}
&\E(\bar{J}_{1}(t))\leq\E\int_{0}^{t} C(1+|\bar z(s)|^{2}+\bar{z}(s-\tau)|^{2})ds\\
&\leq C + C\int_{0}^{t}\sup_{0\leq u\leq s}\E(|z(u)|^{2})ds.
\end{split}
\end{equation}
Using the definition of $z(s)$ and $\bar{z}(s),$ we have
\begin{equation}\label{H402}
\begin{split}
&\E(\bar{J}_{2}(t)) = 2\E\int_{0}^{t}\langle\int_{[\frac{s}{h}]h}^{s}f_{h}(\bar{z}(r),\bar{z}(r-\tau))dr,f_{h}(\bar{z}(s),\bar{z}(s-\tau)\rangle ds\\
&+2\E\int_{0}^{t}\langle\int_{[\frac{s}{h}]h}^{s}g_{h}(\bar{z}(r-),\bar{z}((r-\tau)-),u)\tilde{N}(du,dr),f_{h}(\bar{z}(s),\bar{z}(s-\tau)\rangle ds\\
&=2\E\int_{0}^{t}\langle\int_{[\frac{s}{h}]h}^{s}f_{h}(\bar{z}(r),\bar{z}(r-\tau))dr,f_{h}(\bar{z}(s),\bar{z}(s-\tau)\rangle ds\\
&+2\E\int_{0}^{t}\langle\int_{[\frac{s}{h}]h}^{s}g_{h}(\bar{z}(r-),\bar{z}((r-\tau)-),u)\tilde{N}(du,dr)\big|\mathcal{F}_{[\frac{s}{h}]h},f_{h}(\bar{z}(s),\bar{z}(s-\tau)\rangle ds\\
&=2\E\int_{0}^{t}\langle\int_{[\frac{s}{h}]h}^{s}f_{h}(\bar{z}(r),\bar{z}(r-\tau))dr,f_{h}(\bar{z}(s),\bar{z}(s-\tau)\rangle ds\\
&\leq cth^{1-2\alpha}\leq C.
\end{split}
\end{equation}

By using \eqref{Taylor}, we have 
\begin{equation}\label{H403}
\begin{split}
&\E(\bar{J}_{3}(t))\leq C\E\int_{0}^{t}(1+|\bar z(s)|^{2}+|\bar{z}(s-\tau)|^{2})ds\\
&\leq C + C\int_{0}^{t}\sup_{0\leq u\leq s}\E(|z(u)|^{2})ds.
\end{split}
\end{equation}
By taking the expectation of $\bar{J}_{4}(t),$ we can write that $\E(\bar{J}_{4}(t)) = 0.$ Now recall \eqref{H304}, we have
\begin{equation*}
\begin{split}
&\sup_{0\leq t \leq T}\E(|z(t)|^{2}) \le C+ C\big(\sup_{0\leq t \leq T}\E(|z(t)-D(\bar{z}(t-\tau))|^{2})\big)\\
&\leq C+C\int_0^T\sup_{0\le u\le s}\E(|z(u)|^2)ds.
\end{split}
\end{equation*}
The required result \eqref{H404} will follow an application of the Gronwall inequality.
Now, letting $p=4$ and using the It\^o formula, we have
 \begin{equation}\label{F1}
 \begin{split}
 &|z(t)-G(\bar{z}(t-\tau))|^{p} = |\xi(0)-G(\xi(-\tau))|^{p}\\
 &+p\int_{0}^{t}|z(s)-G(\bar{z}(s-\tau))|^{p-2}\langle z(s)-G(\bar{z}(s-\tau)),f_{h}(\bar{z}(s),\bar{z}(s-\tau))\rangle ds\\
 &+\bigg(\int_{0}^{t}\int_{U}|z(s)-G(\bar{z}(s-\tau))+g(\bar{z}(s),\bar{z}(s-\tau),u)|^{p}-|z(s)-G(\bar{z}(s-\tau))|^{p}\\
 &-p|z(s)-G(\bar{z}(s-\tau))|^{p-2}\langle z(s)-G(\bar{z}(s-\tau)),g(\bar{z}(s),\bar{z}(s-\tau),u)\rangle\bigg) \lambda(du)ds\\
 &+\int_{0}^{t}\int_{U}p|z(s-)-G(\bar{z}((s-\tau)-))|^{p-2}\langle z(s-)-G(\bar{z}((s-\tau)-)),\\
 &\quad \quad \quad \quad \quad \quad \quad \quad \quad \quad \quad \quad \quad \quad \quad \quad \quad \quad \quad \quad g(\bar{z}(s-),\bar{z}((s-\tau)-),u)\rangle\tilde{N}(du,ds)\\
 &=: F_{1}(t)+F_{2}(t)+F_{3}(t)+F_{4}(t).
 \end{split}
 \end{equation}
Noting that
 \begin{equation*}
 \begin{split}
 F_{2}(t)&=p\int_{0}^{t}|z(s)-G(\bar{z}(s-\tau))|^{p-2}\langle \bar{z}(s)-G(\bar{z}(s-\tau)),f_{h}(\bar{z}(s),\bar{z}(s-\tau))\rangle ds\\
 &+p\int_{0}^{t}|z(s)-G(\bar{z}(s-\tau))|^{p-2}\langle z(s)-\bar{z}(s),f_{h}(\bar{z}(s),\bar{z}(s-\tau))\rangle ds,
 \end{split}
 \end{equation*}
 using assumption {\bf(B1)},  \eqref{H404} and \eqref{ErrorWithDriftJump}, we arrive at
 \begin{equation}\label{F2}
 \begin{split}
 &\E(\sup_{0\leq t\leq T}F_{2}(t))\leq C\E\int_{0}^{T}|z(s)-G(\bar{z}(s-\tau))|^{p-2}(1+|\bar z(s)|^{2}+|\bar z(s-\tau)|^{2})ds\\
 & +\E(\sup_{0\leq t\leq T}\int_{0}^{T}|z(s)-G(\bar{z}(s-\tau))|^{p-2}|(z(s)-\bar{z}(s))f_{h}(\bar{z}(s),\bar{z}(s-\tau))|)ds\\
 &\leq C\E\int_{0}^{T}|z(s)-G(\bar{z}(s-\tau))|^{p-2}(1+|\bar z(s)|^{2}+|\bar z(s-\tau)|^{2})ds\\
 &+ \frac{1}{4}\E(\sup_{0\leq t\leq T}|z(t)-G(\bar{z}(t-\tau))|^{p})+\int_{0}^{T}\E(\sup_{0\leq u\leq s}|z(u)-\bar{z}(u)|^{p/2}|f_{h}(\bar{z}(u),\bar{z}(u-\tau))|^{p/2})ds\\
 &\leq C\int_{0}^{T}\E\big(\sup_{0\leq u\leq s}|z(u)|^{p}\big)ds+ \frac{1}{4}\E(\sup_{0\leq t\leq T}|z(t)-G(\bar{z}(t-\tau))|^{p})+C.
 \end{split}
 \end{equation}
By \eqref{Taylor}, {\bf(B1)} and {\bf(B2)}, we obtain
 \begin{equation}\label{F3}
 \begin{split}
 &\E(\sup_{0\leq t\leq T}F_{3}(t))=\E \int_{0}^{T}\int_{U}|z(s)-G(\bar{z}(s-\tau))+g(\bar{z}(s),\bar{z}(s-\tau),u)|^{p}-|z(s)-G(\bar{z}(s-\tau))|^{p}\\
 & -p|z(s)-G(\bar{z}(s-\tau))|^{p-2}\langle z(s)-\bar{z}(x(s-\tau)),g(\bar{z}(s),\bar{z}(s-\tau),u)\rangle \lambda(du)ds\\
 &\leq \E\int_{0}^{T}\int_{U}|z(s)-G(\bar{z}(s-\tau))|^{p-2}|g(\bar{z}(s),\bar{z}(s-\tau),u)|^{2}+|g(\bar{z}(s),\bar{z}(s-\tau),u)|^{p}\lambda(du)ds\\
 &\leq \E\int_{0}^{T}|z(s)-G(\bar{z}(s-\tau))|^{p-2}(1+|\bar{z}(s)|^{2}+|\bar{z}(s-\tau)|^{2})+(1+|\bar{z}(s)|^{p}+|\bar{z}(s-\tau)|^{p})ds\\
 &\leq \frac{1}{4}\E(\sup_{0\leq t\leq T}|z(t)-G(\bar{z}(t-\tau))|^{p})+\E\int_{0}^{T}1+|\bar{z}(s)|^{p}+|\bar{z}(s-\tau)|^{p}ds\\
 &\leq \frac{1}{4}\E(\sup_{0\leq t\leq T}|z(t)-G(\bar{z}(t-\tau))|^{p})+C+  C\int_{0}^{T}\E\big(\sup_{0\leq u\leq s}|z(u)|^{p}\big)ds.
 \end{split}
 \end{equation}
Using  Lemma \ref{BJ} and Young's inequality, we derive that
\begin{equation}\label{F4}
\begin{split}
&\E(\sup_{0\leq t\leq T}F_{4}(t))\\
&\leq C\E\bigg(\int_{0}^{T}\int_{U}|z(s)-G(\bar{z}(s-\tau))|^{2p-2}\K|g(\bar{z}(s),\bar{z}(s-\tau),u)|^{2} \lambda(du)ds\bigg)^{1/2}\\
&\leq C\E[(\sup_{0\leq t\leq T}|z(t)-G(\bar{z}(t-\tau))|^{p-1})\K\big(\int_{0}^{T}\int_{U}|g(\bar z(s),\bar{z}(s-\tau),u)|^{2}\lambda(du)ds\big)^{1/2}]\\
&\leq C\E\int_{0}^{T}|z(s)-G(\bar{z}(s-\tau))|^{p-2}(1+|\bar z(s)|^{2}+|\bar{z}(s-\tau)|^{2})ds\\
&\leq \frac{1}{4}\E(\sup_{0\leq t\leq T}|z(t)-G(\bar{z}(t-\tau))|^{p})+C\E\big(\int_{0}^{T}1+|z(s)|^{p}+|\bar{z}(s)|^{p}+|\bar{z}(s-\tau)|^{p}ds\big)\\
&\leq \frac{1}{4}\E(\sup_{0\leq t\leq T}|z(t)-G(\bar{z}(t-\tau))|^{p})+C+ C\int_{0}^{T}\E\big(\sup_{0\leq u\leq s}|z(u)|^{p}\big)ds.
\end{split}
 \end{equation}
 Now substituting  \eqref{F2}, \eqref{F3} and \eqref{F4} into \eqref{F1}, we obtain that 
 \begin{equation*}
 \E(\sup_{0\leq t\leq T}|z(t)-G(\bar{z}(t-\tau))|^{p})\leq C+ C\int_{0}^{T}\E\big(\sup_{0\leq u\leq s}|z(u)|^{p}\big)ds. 
 \end{equation*}
 By applying \cite[Lemma 6.4.4]{Mao} (also see \eqref{Lemma301}), we  derive that 
 \begin{equation*}
 \E\big(\sup_{0\leq t\leq T}|z(t)|^{p}\big)\leq C+C\int_{0}^{T}\E\big(\sup_{0\leq u\leq s}|z(u)|^{p}\big)ds.
 \end{equation*}
 By the Gronwall inequality, we have for $p=4,$
 \begin{equation*}
 \E\big(\sup_{0\leq t\leq T}|z(t)|^{p}\big)\leq C.
 \end{equation*}
 By repeating the same procedure, the desired result \eqref{UpperBuondLemma3Jump} can be obtained by induction.
 \hfill $\Box$

\subsection{Proof of the Main Results}

In this subsection, we shall prove our main results.

\noindent
\textbf{Proof of Theorem \ref{ConvergenceRateJump}:} 
Recall that
\begin{equation}\label{eq40}
\begin{split}
|x(t)-z(t)|^{p}&\le \kappa|x(t-\tau)-\bar{z}(t-\tau)|^p\\
&\leq \kappa|x(t)-\bar{z}(t)|^pI_{[-\tau\leq t\leq 0]}(t)+\kappa|x(t)-\bar{z}(t)|^pI_{[0\leq t\leq T]}(t)\\
&+\frac{1}{(1-\kappa)^{p-1}}|x(t)-G(x(t-\tau))-z(t)+G(\bar{z}(t-\tau))|^{p},
\end{split}
\end{equation} 
which implies that
 \begin{equation}\label{RateMainTheorem02Jump}
 \begin{split}
 \E\big(\sup_{0\leq t\leq T}|x(t)-z(t)|^{p}\big)&\leq \frac{1}{(1-\kappa)^{p}}\E\big(\sup_{0\leq t\leq T}|x(t)-G(x(t-\tau))\\
 &-z(t)+G(\bar{z}(t-\tau))|^{p}\big)+Ch^p.
 \end{split}
 \end{equation}
An application of the It\^o formula yields that 
for any $p\geq 2,$ 
\begin{equation}\label{B0}
\begin{split}
&|x(t)-G(x(s-\tau))-z(t)+G(\bar{z}(t-\tau))|^{p}\\
&= p\int_{0}^{t}|x(s)-G(x(s-\tau))-z(s)+G(\bar{z}(s-\tau))|^{p-2}\\
&\K\langle x(s)-G(z((s-\tau))-z(s)+G(\bar{z}(s-\tau)),\\
&~~~~~~~~~~~~~~f(x(s),x(s-\tau))-f_{h}(\bar{z}(s),\bar{z}(s-\tau))\rangle ds\\
&+\int_{0}^{t}\int_{U}|x(s)-G(x(s-\tau))-z(t)+G(\bar{z}(t-\tau))\\
&~~~~~~~~~~~~~+(g(x(s),x(s-\tau),u)-g(\bar{z}(s),\bar{z}(s-\tau),u))|^{p}\\
&-|x(s)-G(x(s-\tau))-z(t)+G(\bar{z}(t-\tau))|^{p}\\
&~~~~~~~~~~~~~-p|x(s)-G(x(s-\tau))-z(s)+G(\bar{z}(s-\tau))|^{p-2}\\
&\K\langle x(s)-G(x(s-\tau))-z(s)+G(\bar{z}(s-\tau)),\\
&~~~~~~~~~~~~~~g(x(s),x(s-\tau),u)-g(\bar{z}(s),\bar{z}(s-\tau),u)\rangle \lambda(du)ds\\
&+\int_{0}^{t}\int_{U}|x(s)-G(x(s-\tau))-z(s)+G(\bar{z}(s-\tau))|^{p-2}\\
&\K\langle x(s)-G(x(s-\tau))-z(s)+G(\bar{z}(s-\tau)),\\
&~~~~~~~~~~~~~g(x(s-),x((s-\tau)-),u)-g(\bar{z}(s-),\bar{z}((s-\tau)-),u)\rangle\tilde{N}(du,ds)\\
& =: \bar{F}_{1}(t)+ \bar{F}_{2}(t)+ \bar{F}_{3}(t)
\end{split}
\end{equation}
Using the similar way  as in the \eqref{I1}, \eqref{I2} and \eqref{I3} , we derive that
\begin{equation}\label{B1}
\begin{split}
&\E(\sup_{0\leq t\leq T}\bar{F}_{1}(t))\leq C\int_{0}^{T}\E\big(\sup_{0\leq u\leq s}|x(u)-z(u)|^{p}\big)ds+C\int_{-\tau}^{0}\E\big(\sup_{-\tau\leq\theta\leq 0}|x(\theta)-\bar{z}(\theta)|^{p}\big)ds\\
&+ \frac{1}{2}\E(\sup_{0\leq s\leq T}|x(s)-G(x(s-\tau))-z(s)+G(\bar{z}(s-\tau))|^{p})\\
&+C[(\int_{0}^{T}\E(1+\sup_{0\leq u\leq s}|z(s)|^{l})^{2p}ds)]^{1/2}[\int_{0}^{T}\E(|z(s)-\bar{z}(s)|^{2p})ds]^{1/2}\\
&+C\int_{0}^{T}\E|f(\bar{z}(u),\bar{z}(u-\tau))-f_{h}(\bar{z}(u),\bar{z}(u-\tau))|^{p}ds\\
&\leq C\int_{0}^{T}\E\big(\sup_{0\leq u\leq s}|x(u)-z(u)|^{p}\big)ds +Ch^{p}+Ch^{1/2}+Ch^{\alpha p}\\
&+\frac{1}{2}\E(\sup_{0\leq s\leq T}|x(s)-G(x(s-\tau))-z(s)+G(\bar{z}(s-\tau))|^{p})\\
&\leq C\int_{0}^{t}\E[\sup_{0\leq u\leq s}|x(u)-z(u)|^{p}]ds\\
&+ \frac{1}{2}\E(\sup_{0\leq s\leq t}|x(s)-G(x(s-\tau))-z(s)+G(\bar{z}(s-\tau))|^{p})+Ch^{\gamma},
\end{split}
\end{equation}
where $\gamma = 1/2\wedge \alpha p.$

\noindent
By {\bf(B3)}, {\bf(B4)} and \eqref{Taylor}, we compute
\begin{equation}\label{B2}
\begin{split}\
&\E(\sup_{0\leq t\leq T}\bar{F}_{2}(t))\\
&= \E \int_{0}^{t}\int_{U}|x(s)-G(x(s-\tau))-z(t)+G(\bar{z}(t-\tau))\\
&~~~~~~~~~~~~~+(g(x(s),x(s-\tau),u)-g(\bar{z}(s),\bar{z}(s-\tau),u))|^{p}\\
&-|x(s)-G(x(s-\tau))-z(t)+G(\bar{z}(t-\tau))|^{p}\\
&~~~~~~~~~~~~~-p|x(s)-G(x(s-\tau)x)-z(s)+G(\bar{z}(s-\tau))|^{p-2}\\
&\K\langle x(s)-G(x(s-\tau))-z(s)+G(\bar{z}(s-\tau)),\\
&~~~~~~~~~~~~~~g(x(s),x(s-\tau),u)-g(\bar{z}(s),\bar{z}(s-\tau),u)\rangle \lambda(du)ds\\
&\leq \E\int_{0}^{T}\int_{U}|x(s)-G(x(s-\tau))-z(t)+G(\bar{z}(t-\tau))|^{p-2}\\
&\K|g(x(s),x(s-\tau),u)-g(\bar{z}(s),\bar{z}(s-\tau),u)|^{2}\\
&+|g(x(s),x(s-\tau),u)-g(\bar{z}(s),\bar{z}(s-\tau),u)|^{p}\lambda(du)ds\\
&\leq \E\int_{0}^{T}|x(s)-G(x(s-\tau))-z(s)+G(\bar{z}(s-\tau))|^{p-2}\\
&\K(|x(s)-z(s)|^{2}+|x((s-\tau)-\bar{z}(s-\tau)|^{2})ds\\
&+\int_{0}^{T}|x(s)-z(s)|^{p}+|x(s-\tau)-\bar{z}(s-\tau)|^{p}ds\\
&\leq Ch^{p}+ C\int_{0}^{t}\E[\sup_{0\leq u\leq s}|x(u)-z(u)|^{p}]ds.
\end{split}
\end{equation}
and 
\begin{equation}\label{B3}
\begin{split}
&\E(\sup_{0\leq t\leq T}\bar F_{3}(t))\\
&\leq C\E\bigg(\int_{0}^{T}\int_{U}|x(s)-G(x(s-\tau))-z(s)+G(\bar{z}(s-\tau))|^{2p-2}\\
&\K |g(x(s),x(s-\tau),u)-g(\bar{z}(s),\bar{z}(s-\tau),u)|^{2}\lambda(du)ds\bigg)\\
&\leq  C\E\int_{0}^{T}|x(s)-G(x(s-\tau))-z(s)+G(\bar{z}(s-\tau))|^{p-2}\\
&\K(|x(s)-z(s)|^{2}+|x(s-\tau)-\bar{z}(s-\tau)|^{2})ds\\
&\leq Ch^{p}+C\int_{0}^{t}\E[\sup_{0\leq u\leq s}|x(u)-z(u)|^{p}]ds.
\end{split}
\end{equation}
Now, substituting \eqref{B1}, \eqref{B2} and \eqref{B3} into \eqref{B0}, and using \eqref{RateMainTheorem02Jump}, we have
\begin{equation}\label{RateMainTheorem03Jump}
 \begin{split}
 \E\big(\sup_{0\leq t\leq T}|x(t)-z(t)|^{p}\big)&\leq Ch^{\gamma}+Ch^{p}+C\E\int_{0}^{t}\sup_{0\leq u\leq s}|x(u)-z(u)|^{p}ds\\
 &\leq Ch^{\gamma}+C\int_{0}^{t}\E[\sup_{0\leq u\leq s}|x(u)-z(u)|^{p}]ds.\\
 \end{split}
 \end{equation}
An application of the Grownwall inequality yields the desired result. \hfill
$\Box$

\textbf{Proof of Theorem \ref{ConvergenceJump}:} Noting the fact that all sample paths associated with \eqref{NSDDEJump} are discontinuous,  we define the following stopping times:
For every $R >0,$ stopping times are defined as:
\begin{align}\label{StopppingTimesJump}
\tau_{R}:= \inf\{t\geq 0: |x(t)| > R\}, \quad \rho_{R}:= \inf\{t\geq 0: |z(t)| > R\}. 
\end{align}
The remainder of the proof  is similar to the that of Theorem \ref{Convergence},   we omit details here. \hfill $\Box$

\end{document}